\documentclass[12pt]{article}
\usepackage{epsfig}
\usepackage{rotating}
\title{Le Retour de Pappus}
\author{Richard Evan Schwartz \thanks{\hskip 5 pt 
Supported by  N.S.F. Grant DMS-2102802, a Simons Sabbatical
Fellowship, and a Mercator Fellowship.}}

\newtheorem{theorem}{Theorem}[section]

\newtheorem{lemma}[theorem]{Lemma}

\def\startproof{{\bf {\medskip}{\noindent}Proof: }}

\def\endproof{$\spadesuit$  \newline}

\def\H{\mbox{\boldmath{$H$}}}% 
\def\P{\mbox{\boldmath{$P$}}}% 
\def\R{\mbox{\boldmath{$R$}}}% 
\def\Z{\mbox{\boldmath{$Z$}}}% 

\begin{document}

\maketitle

\begin{abstract}
  In [{\bf S0\/}]
    we explained how the iteration of Pappus's
  Theorem gives rise to a $2$-parameter family
  of representations of the
  modular group into the
  group of projective automorphisms.
    In this paper we realize
    these representations as isometry groups
      of patterns of geodesics
  in the symmetric space $X=SL_3(\R)/SO(3)$.  The
  patterns have the same asymptotic structure as the geodesics in the
  Farey triangulation, so  our construction gives a
  $2$ parameter family of deformations of the
  Farey triangulation inside $X$.  We will also describe
  a bending phenomenon associated to these patterns.
      \end{abstract}

  \section{Introduction}

Pappus's Theorem says that if
$A=(A_1,A_2,A_3)$ and $B=(B_1,B_2,B_2)$
consist of collinear points, then so does
$C=(C_1,C_2,C_3)$.  Starting with the
pair $(A,B)$ you produce the pairs $(A,C)$ and $(C,B)$.

 \begin{center}
\resizebox{!}{1.5in}{\includegraphics{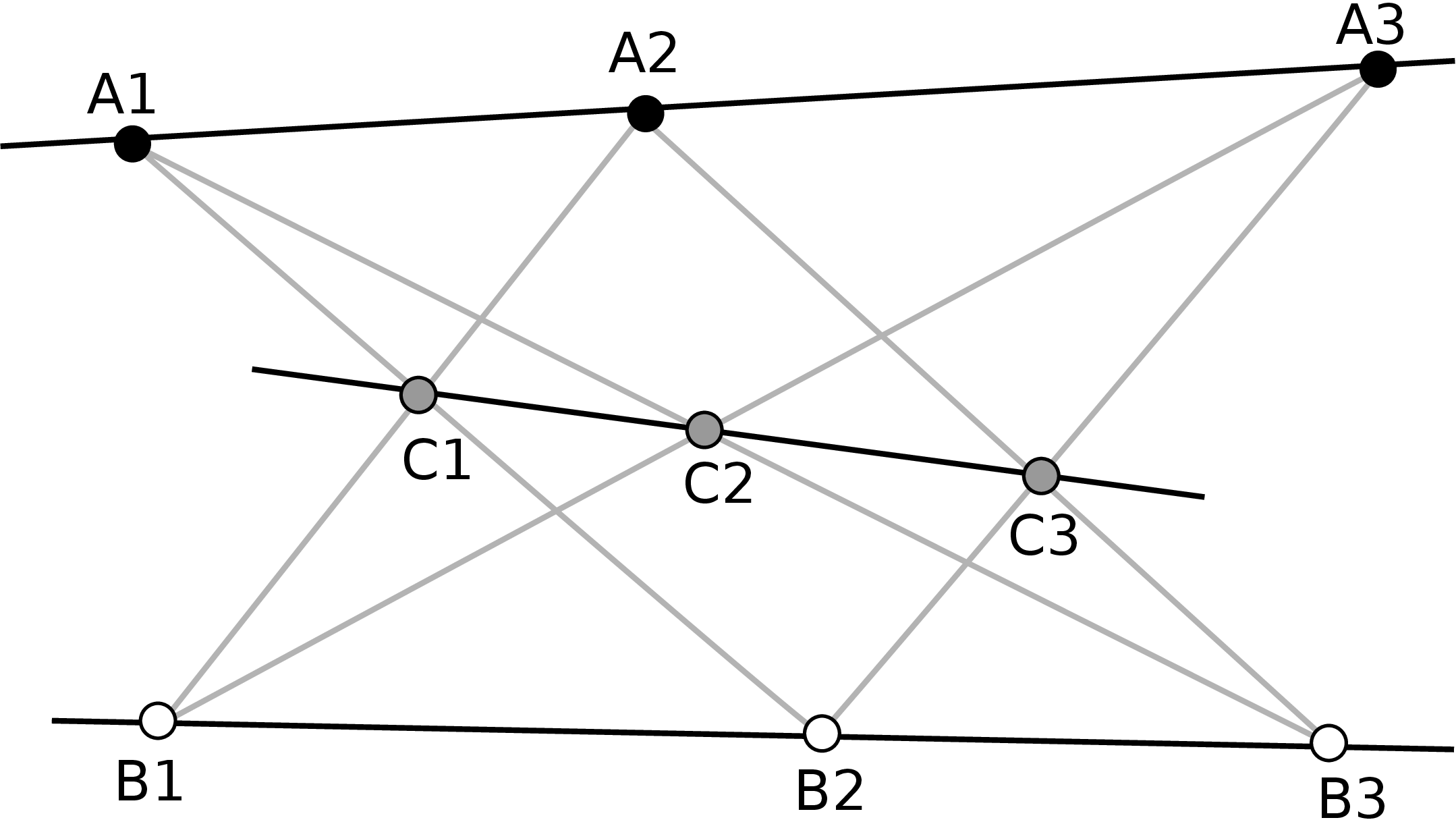}}
\newline
{\bf Figure 1.1:\/} Pappus's Theorem
\end{center}

In my 1993 paper [{\bf S0\/}], I iterate the
above construction and then use the result to define a
$2$ dimensional moduli space of inequivalent
representations of the modular group
into the group of projective symmetries of the
flag variety $\cal P$ over the projective plane.
Here $\cal P$ is the space
of pairs $(p,\ell)$ where $p$ is a one dimensional
subspace in $\R^3$, and $\ell$ is a $2$-dimensional
subspace in $\R^3$, and $p \subset \ell$.
Figure 1.2, a concrete example, hints at the
limit set of one of these {\it Pappus modular groups\/}.

 \begin{center}
\resizebox{!}{4.1in}{\includegraphics{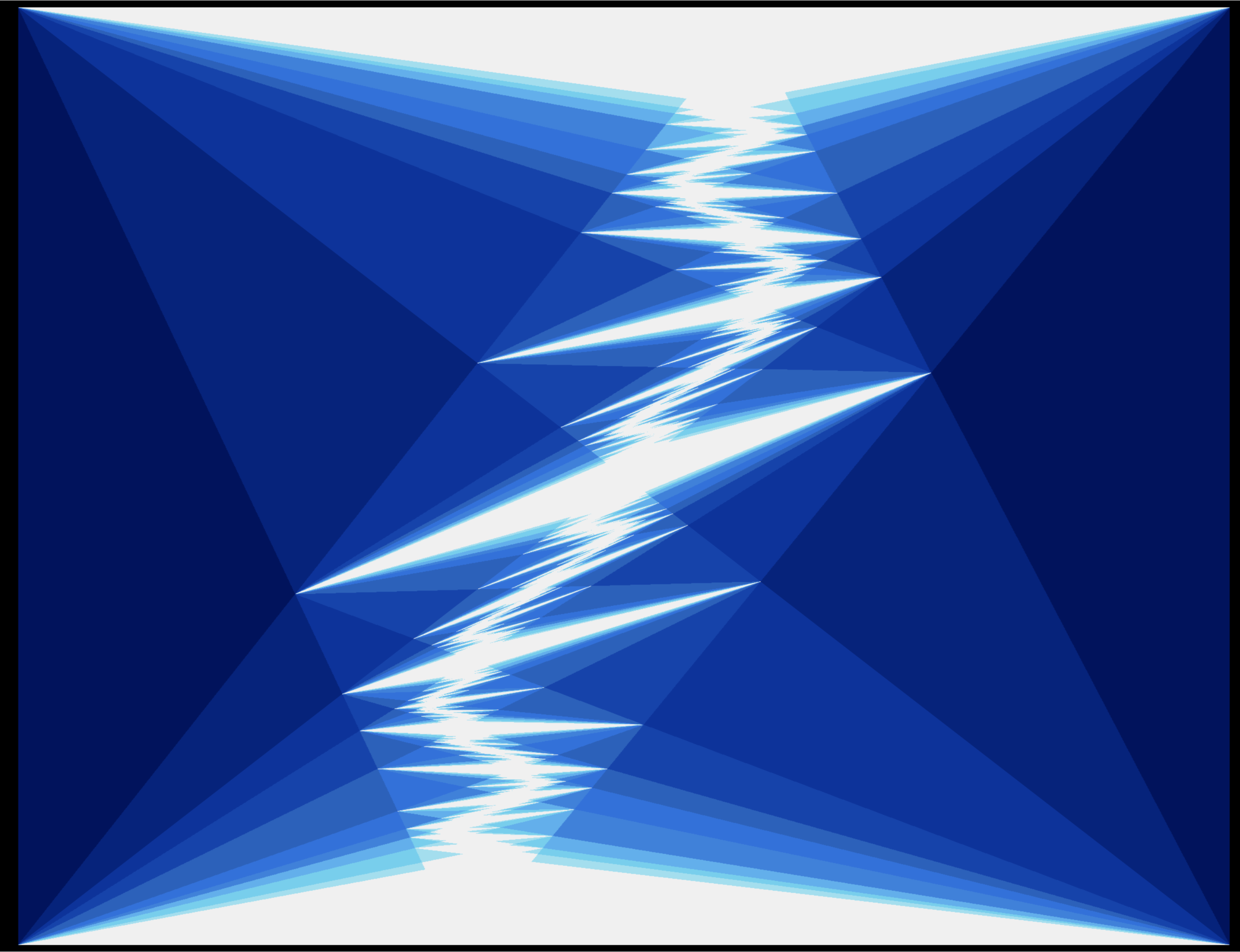}}
\newline
{\bf Figure 1.2\/} A hint of the limit set
\end{center}

The paper
[{\bf S0\/}] is
a precursor of many papers
in higher Teichmuller Theory.
One can see many of the basic structures
in [{\bf S0\/}] appearing much more generally
e.g. in [{\bf Bar\/}], [{\bf BCLS\/}], [{\bf GW\/}],  [{\bf KL\/}], and
[{\bf Lab\/}].
  Nowadays, the Pappus modular groups are
  classified as
  {\it relatively Anosov representations
    in the Barbot component\/}.  Compare
  [{\bf BLV\/}], [{\bf KL\/}], and [{\bf Bar\/}].
See also [{\bf G\/}] and  [{\bf Hit\/}], two papers roughly
contemporaneous with (and even a bit earlier than)
[{\bf S0\/}], that inspired a huge amount of
work in higher Teichmuller Theory.

The purpose of this paper is to give a
description of the Pappus modular groups
in terms of the associated symmetric space,
$X=SL_3(\R)/SO(3)$.  The 
description adds
another layer of depth and beauty to the
Pappus modular groups, and also it
is a textbook picture
of how a relatively Anosov representation
acts on the associated symmetric space as
the symmetries of something like a pleated plane.
I have wanted to do this for many years, but only recently
saw the answer, implicit in 
[{\bf S0\/}, Figure 2.4.2],
staring me in the face.
See Figure 4.2 below.

 \begin{center}
\resizebox{!}{2.5in}{\includegraphics{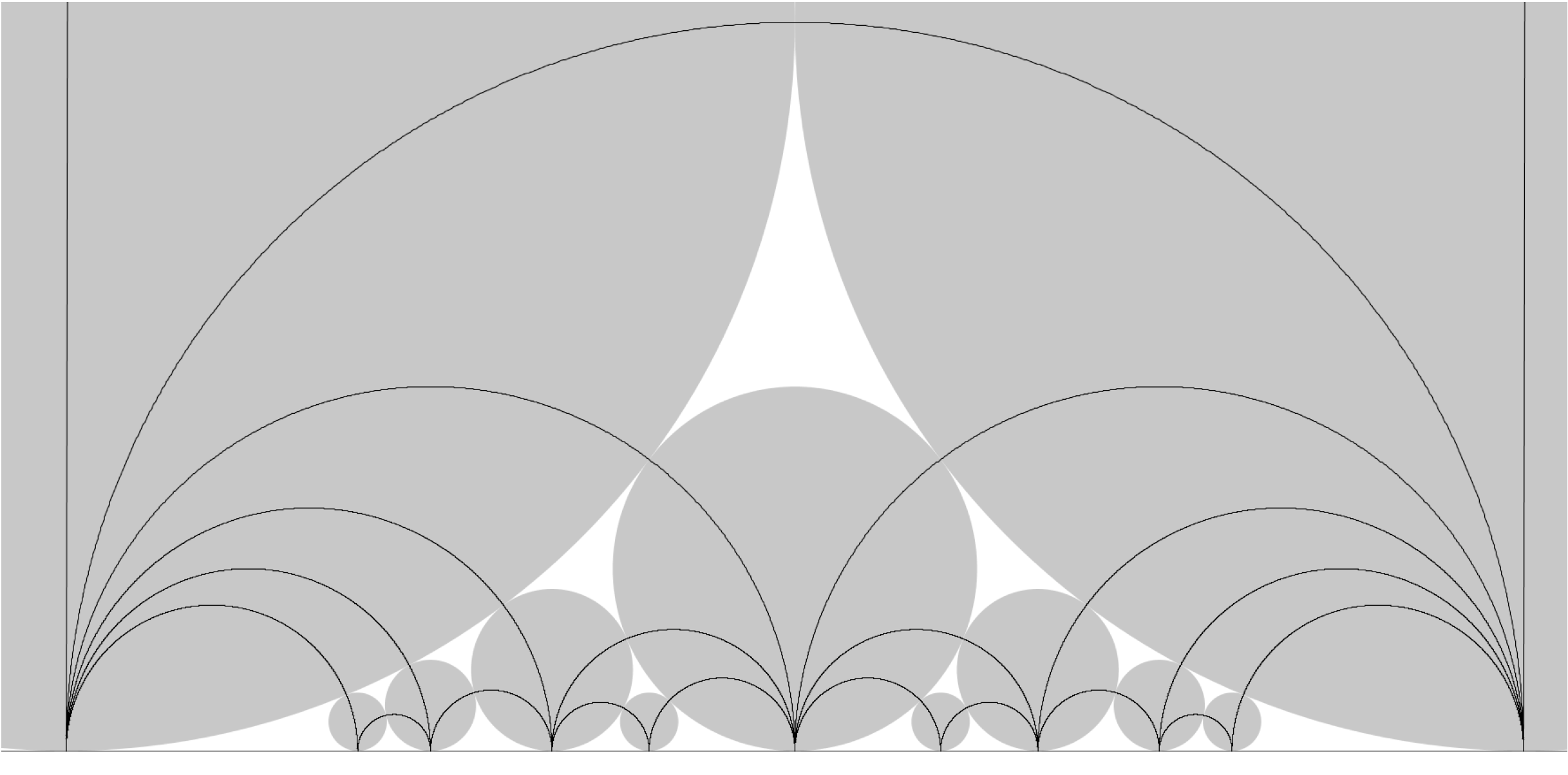}}
\newline
{\bf Figure 1.3:\/} Part of the Farey triangulation and dual horodisk
packing
\end{center}

Working in the upper half-plane model $\H^2$ of the hyperbolic plane,
one forms a pattern of geodesics by connecting each pair of
rational numbers $a/b$ and $c/d$ on the ideal boundary
$\R \cup \infty$ by a geodesic if and only if $|ad-bc|=1$.
Here $\infty$ is interpreted as the fraction $1/0$.  The resulting
pattern of geodesics is the $1$-skeleton of the famous
{\it Farey triangulation\/}, and the
modular group is precisely the group of orientation preserving
symmetries of this thing.

The visual boundary $\partial X$ of $X$ is
$4$-dimensional, and it contains the
$3$-dimensional $\cal P$ in a natural way.
See \S \ref{visb}.
We say that a {\it medial geodesic\/} in $X$ is a geodesic which
limits to points in $\cal P$ at either end.
We say that two medial geodesics are {\it one-end-asymptotic\/}
if, in exactly one direction, they limit on the same flag.
We say that they are {\it non-asymptotic\/} if they have
no ends in common. 
(Note that in a higher rank symmetric space like $X$ one
has distinct geodesics that are asymptotic at both ends.)
\newline
\newline
{\bf Main Result:\/}
Let $F$ denote the union of the geodesics in the
Farey triangulation.  We say that a
{\it Farey pattern\/} in $X$ is an
injective map $\mu: F \to X$ which is an
isometry when restricted to each geodesic
of $F$ and which maps pairs of one-end-asymptotic geodesics
to pairs of one-end-asymptotic geodesics and pairs of non-asymptotic
geodesics to pairs of non-asymptotic geodesics.

\begin{theorem}
  \label{pattern}
  Every Pappus modular group arises as a
  group of isometries of a Farey pattern.
  Generically, the Pappus modular group is the full isometry
  group of the pattern.
    \end{theorem}

\noindent
{\bf Proof Sketch:\/}
A Pappus modular group $\Lambda_{\cal M}$ is the
group of projective symmetries of a certain infinite collection
$\cal M$ of {\it marked boxes\/}, objects which encode instances
of Pappus's Theorem. Each marked box $M \in \cal M$ has associated to it
a pair of flags $\phi_{M,1}$ and $\phi_{M,2}$
and also an order $2$ element $\rho_M \in \Lambda_{\cal M}$.
It turns out that $\rho_M$, when acting on $X$, has a unique
fixed point $p_M$. 
There is a natural $1$-parameter family $f_M$ of medial geodesics
which foliate a flat in $X$ and
limit at their ends to $\phi_{M,1}$ and $\phi_{M,2}$.  It turns out
that within $f_M$
there is a unique medial geodesic $\gamma_M$ in $f_M$
which contains $p_M$.  The assignment $M \to \gamma_M$,
as $M$ ranges over $\cal M$, gives us the Farey pattern.
\newline

\noindent
{\bf The Bending Phenomenon:\/}
Let $\Gamma$ be one of our Farey patterns.
Each geodesic in $\Gamma$ is contained in a unique flat.
Our proof that the geodesics of $\Gamma$ are disjoint
will actually show that the associated flats are disjoint
as well.  Say that a {\it prism\/} is a triple of flats
corresponding to an ideal triangle in the Farey triangulation.
We want to study the geometry of these prisms.

Referring to Figure 1.1, the quantity $|\log(-\chi)|$, where
$\chi$ is the triple product (as in Equation \ref{tripp})
of the three flags
$$(A_2,\overline{A_1A_3}), \hskip 20pt  (B_2,\overline{B_1B_3}),
\hskip 20 pt
(C_2,\overline{C_1C_3})$$
gives rise to what we call the {\it triple invariant\/} of a
Pappus modular group.  Compare [{\bf FG\/}].
The level sets of the triple invariant give $1$-parameter families of
Pappus modular groups.
Figure 5.3 shows a plot on the character variety.

The geometry of the individual associated prisms
only depends on the triple invariant.  Hence,
within such a $1$-parameter level set,
the prisms are all isometric to
each other.  As we vary within the level set,
we have a bending/shearing phenomenon akin to
those found in [{\bf T\/}],  [{\bf P\/}] and [{\bf FG\/}].
The geometry of the prisms does not change, but the geometry
of pairs of adjacent prisms does change.
\newline
\newline
{\bf An Attempt at a Pleated Surface:\/}
The Farey patterns typically do not lie in a
totally geodesic slice of $X$.  Given the bending
phenomenon just discussed, one might wonder
whether there is a natural way to fill in around
the geodesics to get a kind of pleated surface.
In \S \ref{rough} I will explain one way to do it,
using a coning construction.  This gives rise to
a piecewise analytic
disk with totally geodesic pleats.
Some of the pleats are
bi-infinite geodesics and some are geodesic rays
meeting $3$ at a time.  I do not know how to prove
that this disk is embedded, though I guess that it is.

I will also describe how the coning 
operation gives rise to a kind
of filling of the associated prisms, giving
something like a pleated $3$ manifold.
I do not know how to prove that this
thing is embedded and I don't really see the
picture very well.  I suspect that I don't have
the filling picture in focus yet. Anna Wienhard
suggested to me that perhaps the methods in
[{\bf DR\/}] would do a better job.
\newline

The rest of this paper fleshes out the above sketches.
In \S 2 I will explore the geometry of the symmetric space $X$,
explaining it in very elementary terms.
In \S 3 I will exposit large parts
of [{\bf S0\/}].  To avoid repeating what I already did in [{\bf S0\/}], and also
benefitting from years of hindsight, I tried to give a more
geometric and accessible account here.
(See [{\bf BLV\/}] for another exposition of [{\bf S0\/}] closer
to my original one.)
In \S 4 I will put everything together and prove Theorem
\ref{pattern}.  In \S 5 I will discuss the bending phenomenon
and also write down everything else I have figured out about the
structure of these patterns in $X$.
\newline

I  thank Martin Bridgeman,
Bill Goldman,  Joaquin Lejtreger,
Joaquin Lema, Dan Margalit,  Dennis
Sullivan, and Anna Wienhard for various
interesting and helpful conversations.
I thank the referee for a careful reading of the paper.
Finally, I especially
thank Martin and the two Joaquins
for recently rekindling my interest in the subject.

\newpage

  \section{Geometry of the Symmetric Space}

\subsection{Primer on Projective Geometry}
\label{projgeo}

\noindent
{\bf The Projective Plane:\/}
The {\it projective plane\/} $\P$ is the space of $1$-dimensional
subspaces of $\R^3$.   The {\it affine patch\/} in $\P$ is the
subset consisting of subspaces not contained in the $XY$-plane.
The affine patch is essentially a copy of $\R^2$ sitting in $\P$:
The subspace containing $(x,y,1)$ is identified with $(x,y) \in \R^2$.
\newline
\newline
{\bf The Dual Projective Plane:\/}
The {\it dual projective plane\/} $\P^*$ is
the space of $2$-dimensional subspaces of $\R^3$. A
{\it line\/} in $\P$ is the set of all $1$-dimensional
subspaces contained in a fixed $2$-dimensional subspace.
Thus the lines in $\P$ are in canonical bijection with the
points of $\P^*$.
\newline
\newline
{\bf Projective Transformations:\/}
A {\it projective transformation\/} is a map of $\P$ induced by
an invertible linear transformation.  This notion makes sense
because linear transformations permute the $1$-dimensional
subspaces of $\R^3$.   Projective transformations simultaneously
act on $\P^*$.   Projective transformations act as analytic
diffeomorphisms on $\P$ and on $\P^*$.
Projective transformations of $\P$ maps lines to lines.
The group of projective transformations is $8$ dimensional
and acts simply transitively on the space of quadruples
of general position points in $\P$.  Thus there is a unique
projective transformation which takes any given general
position quaduple to another given one.
\newline
\newline
{\bf Dualities:\/}
A {\it duality\/} is a factor-switching
map $\delta: \P \cup \P^* \to \P^* \cup \P$
which maps collinear points to coincident
lines.  The duality is a {\it polarity\/} if
it has order $2$.   A polarity is induced
by a quadratic form:  Each subspace is
mapped to the subspace perpendicular
to it with respect to the quadratic form.
If this quadratic form is positive definite,
the polarity is called {\it elliptic\/}.
The dot product induces what we call
the {\it standard polarity\/}, which we call $\Delta$.
Any elliptic
polarity has the form $T\circ \Delta \circ T^{-1}$
where $T$ is a projective transformation.
\newline
\newline
{\bf The Flag Variety:\/}
The {\it flag variety\/} is the subspace of
$\P \times \P^*$ consisting of pairs
$(p,\ell)$ where $p \in \ell$.  These
elements are called {\it flags\/}.  We
call this space $\cal P$.  The space
$\cal P$ is a $3$-dimensional manifold.
A projective transformation $T$ and a duality $\delta$
act as analytic diffeomorphisms of $\cal P$ as follows:
$$T(p,\ell)=(T(p),T(\ell)), \hskip 30 pt
\delta(p,\ell)=(\delta(\ell),\delta(p)).$$

\subsection{Ellipsoids}

\noindent
{\bf Frames:\/}
We work with ellipsoids of volume $4 \pi/3$, which
are centered at the origin in $\R^3$.   
We call such ellipsoids {\it unit ellipsoids\/}.  The unit ball is
an example.
We say that a {\it frame\/} is a triple of
mutually orthogonal lines through the origin.
The {\it standard ellipsoids\/} are given by
\begin{equation}
  \label{standard}
  \frac{x^2}{a^2} + \frac{y^2}{b^2}+\frac{z^2}{c^2} \leq 1, \hskip 30 pt
  a,b,c>0, \hskip 20 pt abc=1.
\end{equation}
The {\it standard frame\/} 
is the triple of coordinate axes.

We say that the frame $f$ is {\it associated\/} to the
ellipsoid $E$ if there is an isometry $I \in SO(3)$ such
that $I(E)$ is a standard ellipsoid and $I(f)$ is the standard frame.
The unit ball has an $SO(3)$-family of frames associated to it.
The generic ellipsoid,  e.g. a standard ellipsoid with
$a\not = b \not =c$, has a unique associated frame up to
permutation of the lines.  All other ones, e.g. a standard
ellipsoid with $a=b \not = c$, have an $S^1$-family of
associated frames.
\newline
\newline
{\bf Principal Lengths:\/}
We get the {\it principal lengths\/} of an
ellipsoid by intersecting it with the lines
of any associated frame and taking half the lengths.
Thus, the principal
lengths of the ellipsoid in Equation
\ref{standard} are $a,b,c$.  These principal
lengths are sometimes called {\it singular values\/},
in connection with the singular value decomposition
of a symmetric matrix.

If the principal lengths of an ellipsoid $E$
are $a,b,c$ we define
\begin{equation}
  \lambda(E)=\|V\|=\sqrt{V \cdot V}, \hskip 30 pt
  V=(\log a,\log b,\log c).
\end{equation}
Thus $\lambda(E) \geq 0$, with equality if and only if
$E$ is the unit ball.
\newline
\newline
{\bf Matrix Representation:\/}
Every unit ellipsoid has associated to it a
unit determinant positive definite
symmetric matrix $S$ in the following
way.   The ellipsoid $E_S$ is the set of vectors $v$
such that $S(v) \cdot v \leq 1$.    We say that $S$
{\it represents\/} $E_S$.  The identity matrix represents
the unit ball.  The principal lengths associated to an
ellipsoid are the reciprocals of the square-roots of the
eigenvalues of the
representing matrix.
\newline
\newline
{\bf Group Action:\/}
$SL_3(\R)$ acts on the space $X$ of unit ellipsoids
in the obvious way. If $T \in SL_3(\R)$ and $E \in X$, then
we have $T(E) \in X$ as well.   If the matrix $S$ represents
$E$ then the matrix
$ (T^{-1})^t S T^{-1}$
represents $T(E)$.
Here $(\cdot)^t$ is the transpose.  To see this, we note that
$$(T^{-1})^t S T^{-1}(T(v)) \cdot T(v)=(T^{-1})^t(S(v)) \cdot
T(v)=S(v) \cdot v.$$

\subsection{The Symmetric Space}
\label{symm}

Here we discuss the geometry of the space $X$ of unit
ellispoids
centered at the origin.  We first mention
several other interpretations.  First, we can think of $X$ as the
space of positive definite symmetric matrices of determinant one.
As discussed above, each such matrix defines the ellipsoid which
is its unit ball. 
Second, we can think of $X$ as $SL_3(\R)/SO(3)$.  To see the connection, note that
$SL_3(\R)$ acts transitively on $X$
and the
stabilizer of the unit ball, a distinguished origin in $X$,  is exactly $SO(3)$.
\newline
\newline
{\bf The Metric:\/}
The metric $d$ on $X$ can be described as follows.
\begin{itemize}
\item $d(ME_1,ME_2)=d(E_1,E_2)$ for any $M \in SL_3(\R)$
and $E_1,E_2 \in X$.
\item If $E_1$ is the unit ball then $d(E_1,E_2)=\lambda(E_2)$.
\end{itemize}
By construction, $SL_3(\R)$ acts isometrically on $X$.
It turns out that $d$ is induced by a Riemannian metric
of non-positive sectional curvature.
\newline
\newline
{\bf Geodesics through the Origin:\/}
Let $E \in X$ be any
point other than the unit ball.
Let $L_1,L_2,L_3$ be the lines in an
associated frame, and let
$\lambda_1,\lambda_2,\lambda_3$ be the
corresponding principal lengths.
We write $\ell_j=\log(\lambda_j)$.  We let
$E(t)$ denote the geodesic with the same
associated frame and principal lengths
$\{\exp(t \ell_j)\}$.  By construction,
$E(1)=E$.  The collection $\{E(t)|\ t \in \R\}$ is
a geodesic in $X$ through the origin.
The origin is given by $E(0)$.
In terms of the matrix representations, there is a
$1$-parameter subgroup associated to the
geodesic defined by the unit ball $E_0$ and $E$, and this
subgroup is conjugate in $SO(3)$ to a
$1$-parameter subgroup of positive diagonal matrices.
Every other geodesic in $X$ is carried to a geodesic
through the origin by an element of $SL_3(\R)$.
\newline
\newline
{\bf Flats:\/}
A {\it flat\/} is a totally geodesic copy of $\R^2$ embedded
in $X$.    The {\it standard flat\/} is
the union of the points corresponding to the standard
ellipsoids.
The set of matrices representing the points in the standard
flat is exactly the subgroup of positive diagonal matrices.  So,
the standard flat is the orbit of the origin under
the subgroup of positive diagonal matrices.
The standard flat has $3$ {\it singular geodesics\/} through the origin,
the ones represented by matrices having repeated eigenvalues.
These $3$ singular geodesics divide the standard flat into
$6$ sectors which are called  {\it Weyl chambers\/}.  See Figure 2.1 below.
We have this Weyl chamber structure at
each point of each flat.
\newline
\newline
{\bf The Basic Involution:\/}
Every elliptic polarity defines an isometry of $X$ with a
unique fixed point.   Since the elliptic polarities are all
conjugate in $SL_3(\R)$ we just have to understand this
for the elliptic polarity $\Delta$ induced by the dot product.
For any ellipsoid $E$, the ellipsoid $\Delta(E)$ is the
ellipsoid which has the same associated frame as $E$ and
reciprocal principal lengths.  The matrix representing
$\Delta(E)$ is the inverse of the matrix representing $E$.
We can see quite readily that $\Delta$ fixes the origin
in $X$ and reverses all the geodesics through the origin.
Since $X$ is a symmetric space, the geodesic-reversal property
implies that $\Delta$ acts as an isometry.
To be sure, let me give a self-contained proof.

\begin{lemma}
  $\Delta$ is an isometry of $X$.
\end{lemma}

\startproof
Let $Y^*=(Y^{-1})^t$ for any
$Y \in SL_3(\R)$.  The map $Y \to Y^*$ is an order $2$
automorphism of $SL_3(\R)$.

Let $E_0$ be the unit ball, as above.
Two typical points in $E_1,E_2 \in X$ are given by
$M_1(E_0)$ and $M_2(E_0)$ where $M_1,M_2 \in SL_3(\R)$.
The positive definite symmetric matrix representing
$E_k$ is $S_k=M_k^*M_k^{-1}$.    The matrix
representing $\Delta(E_k)$ is
$$S_k^{-1}=M_k(M_k^*)^{-1}=M_kM_k^t.$$
From this we recognize that $\Delta(E_k)=M_k^*(E_0)$.

From the $SL_3(\R)$-symmetry of the metric, we have
$$d(E_1,E_2)=d(E_0,M_1^{-1}(E_2))=d(E_0,M_1^{-1}(M_2(E_0)))=\lambda(M_1^{-1}M_2(E_0)).$$
Likewise
$$d(\Delta(E_1),\Delta(E_2))=\lambda(M_1^tM_2^*(E_0)).$$
It is convenient to set
$$M=M_1^{-1}M_2, \hskip 30 pt S=M^*M^{-1}.$$
With this notation,
$$d(E_1,E_2)=\lambda(M(E_0)), \hskip 30 pt
d(\Delta(E_1),\Delta(E_2))=
\lambda(M^*(E_0)).$$
The matrix $S$ represents $M(E_0)$.  Because $(*)$
is an automorphism, $S^*$ represents $M^*(E_0)$.
But the eigenvalues of $S$ and the eigenvalues of $S^*$
are reciprocals of each other.  Hence
$\lambda(M(E_0))=\lambda(M^*(E_0))$.
\endproof

\subsection{The Visual Boundary}
\label{visb}

The {\it visual boundary\/}  of $X$ is the space
$\partial X$
of geodesic rays through the origin.   Topologically, the
visual boundary of $X$ is $S^4$, the $4$-dimensional
sphere.   We identify $3$ special subsets of the visual
boundary:
\begin{itemize}
\item The rays corresponding to ellipsoids having principal lengths
  $$\lambda^2,1/\lambda,1/\lambda,$$
  with $\lambda  \geq 1$.
  As $\lambda \to \infty$, the corresponding ellipsoids
  look like uncooked spaghetti and they converge to
  a line through the origin.
   Thus, the corresponding subset of
 $\partial X$ is a copy of  the projective plane $\P$.
\item The rays corresponding to ellipsoids having principal lengths
  $$\lambda,\lambda,1/\lambda^2,$$
  with $\lambda \geq 1$.
  As $\lambda \to \infty$, 
  such ellipsoids look like pancakes and
    they converge to a unique plane through the origin.
 Thus, the corresponding subset of
 $\partial X$ is a copy of  $\P^*$.
 \item The rays corresponding to ellipsoids having principal lengths
  $$\lambda,1,1/\lambda,$$
  with $\lambda   \geq 1$.
  As $\lambda \to \infty$ the ellipsoids in the same ray
  look like popsicle sticks and their limit defines a unique flag:
  The ellipsoids
  accumulate to a unique line through the origin, one of the
  lines in the associated frame, and their two largest
  principle directions pick out a plane, one that is spanned
  by two of the lines in the associated frame.
   Thus, the corresponding subset of
 $\partial X$ is a copy of  $\cal P$.
 \end{itemize}

 Two general geodesic rays in $X$ are
 {\it fellow travelers\/} if there is a uniform bound on
 the distance between corresponding points on the rays.
 
 \begin{lemma}
   \label{asymptotic}
   Every geodesic ray in $X$ is the fellow traveler of a geodesic ray
   through
   the origin.
 \end{lemma}

 \startproof
 This is well known.
 Let $\gamma$ be such a ray.  We have $\gamma=T(\gamma_0)$
 where $T \in SL_3(\R)$ and $\gamma_0$ is a ray through the origin.
 By symmetry we can assume that $\gamma_0$ is a ray in the standard flat.

 There are $3$ numbers $\ell_1,\ell_2,\ell_3$ with $\ell_1+\ell_2+\ell_3=0$ such that
 the principal lengths of $E_t \in \gamma_0$ is
 $\exp(t \ell_j)$ for $j=1,2,3$.   We order these numbers so that
 $\ell_1 \geq \ell_2 \geq \ell_3$.  Here $t \geq 0$.
 We associate to $\gamma_0$ the standard frame $f_0$.

 Let $L_1,L_2,L_3$ be the lines $T(f_0)$.  These lines
 are not necessarily mutually perpendicular.
 Let $L_1^*=L_1$.  Let $L_2^*$ be the line through the origin
 perpendicular to
 $L_1$ inside the span of $L_1,L_2$.   Let $L_3^*$ be the
 line through the origin mutally perpendicular to
 $L_1^*$ and $L_2^*$.    Let
 $E^*_t$ denote the ellipsoid whose associated
 frame is $(L_1^*,L_2^*,L_3^*)$ and whose
 principal lengths are
 $\exp(t\ell_1), \exp(t\ell_2),\exp(t\ell_3)$.
 The set $\gamma_*=\{E^*_t|\ t \in \R_+\}$ is a
 geodesic ray in $X$ through the origin.
 We want to see that $\gamma$ and $\gamma_*$
 are fellow travelers.
 
We can further normalize by an element of $SO(3)$ so that
$\gamma_*=\gamma_0$.  This is to say that
$L_1^*,L_2^*,L_3^*$ are the coordinate axes.
Now we want to see that $\gamma_0$ and $\gamma$ are
fellow travelers.  Consider the diagonal matrix $M_t$ with
diagonal entries $\exp(-t\ell_j)$ for $j=1,2,3$.
By construction $M_t(E_t^*)$ is the unit ball.

Consider $M_t(E_t)$.  This ellipsoid has volume $4 \pi/3$
and intersects the $X$-axis in a segment of length $2$.
We claim that $M_t(E_t)$ intersects the $XY$-plane in
a set of area $A_t$ that is uniformly bounded away
from $0$ and $\infty$.  Assuming this claim, $M_t(E_t)$
is uniformly close to the unit ball in shape.   Hence
$d(M_t(E_t),M_t(E_t^*))=d(E_t,E_t^*)$ is uniformly bounded.

For the claim, note that $E_t$ intersects $L_1$ and $L_2$ respectively
in segments of length $2\exp(t\ell_1)$ and
$2\exp(t\ell_2)$.  The first of these numbers is the largest.  Also,
there is a fixed angle between
$L_1$ and $L_2$.  Given this situation, we see that
$E_t$ intersects $L_2^*$ in a segment of length
$2A_t \exp(t\ell_2)$ where $A_t \leq 1$ is uniformly bounded away from $0$.
Hence $M_t(E_t)$ intersects the $XY$-plane in a set
of area $A_t$ that is uniformly bounded away from $0$ and $\infty$.
\endproof

Two distinct geodesic rays through the origin are not fellow
travelers.
Thus, Lemma \ref{asymptotic} gives a way to assign a point
in $\partial X$ to any geodesic ray.  This, in turn, gives
a way to extend the action of $SL_3(\R)$ to $\partial X$.
This action is compatible with the action of $SL_3(\R)$ on
each of the sets $\P$ and $\P^*$ and $\cal P$.
The action of a duality on $X$ also extends
to $\partial X$ and is compatible with the
action we already have on $\P \cup \P^*$ and $\cal P$.
\newline
\newline
{\bf Remark:\/}
Consider a {\it medial ray\/}, a geodesic
ray which is asymptotic to a flag in $\cal P$.  The 
ellipsoids look like popsicle sticks far out on the ray.
If we apply a linear transformation then,
far out on the new ray, the ellipsoids
again look like popsicle sticks.  This defines the
extension more concretely in the  special case
of interest to us.

\subsection{Flats and the Visual Boundary}

Now we explain how  flats interact with the visual
boundary.  Consider the standard flat.
Each of the $3$ singular geodesics through the origin limits
in one direction to one of the coordinate axes and
in the other direction to the complementary
coordinate plane.  Thus, the flat has $3$ limit
points in $\P$ and $3$ in $\P^*$.  The $3$
geodesics that lie midway between the singular geodesics,
all of which are medial geodesics,
limit to points in $\cal P$.  Thus the flat
also has $6$ limit points in $\cal P$.
With respect to the circular order on the rays through
the origin in the flat, the flag points interlace with
the points in
$\P$ and $\P^*$, which are themselves interlaced.
Figure 2.1 shows this structure.
The limit points in $\P$ and $\P^*$ and $\cal P$
are respectively denoted by black, white, and grey arrows.

 \begin{center}
\resizebox{!}{2.2in}{\includegraphics{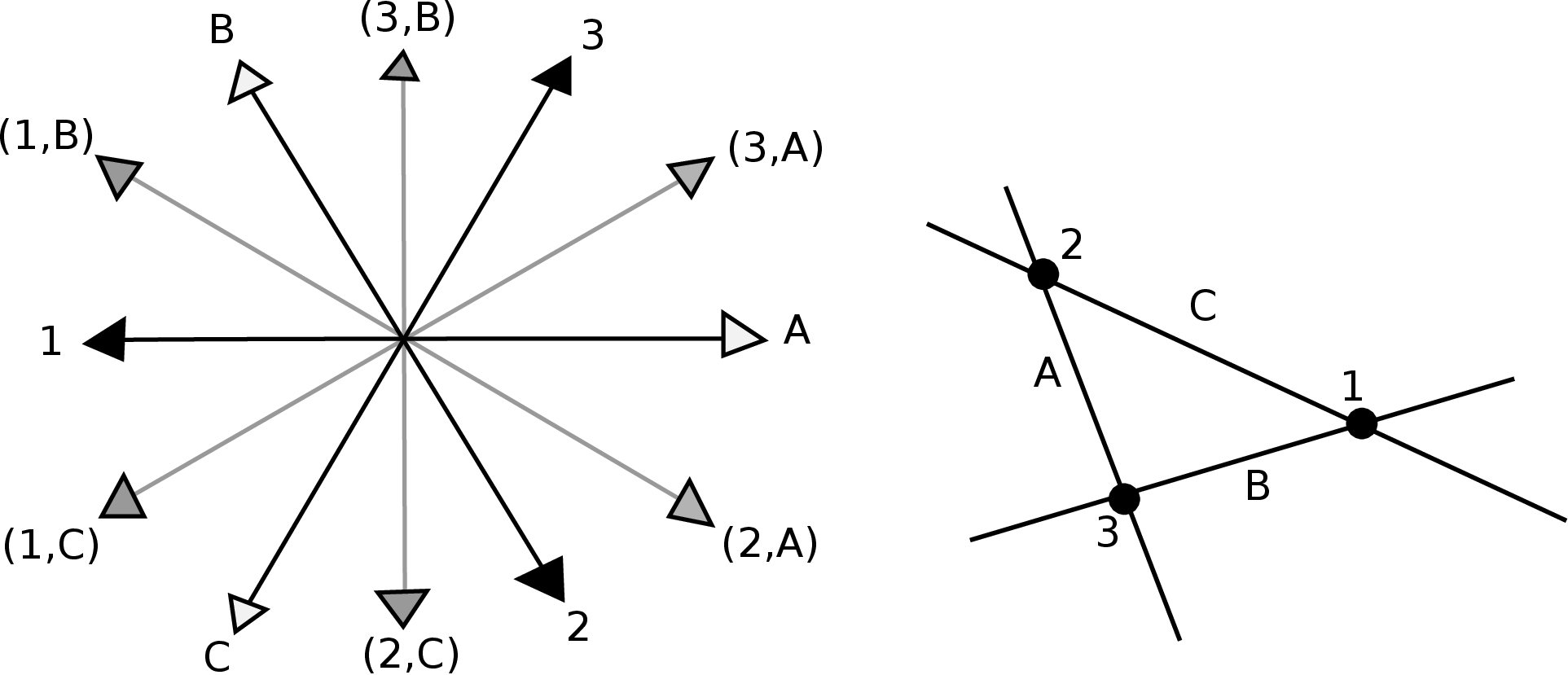}}
\newline
{\bf Figure 2.1:\/}  The picture in the flat (left) and in $\P$
(right).
\end{center}

This union of
$3$ points, $3$ lines, and $6$ flats corresponds
naturally to a triangle in $\P$, and likewise a
triangle in $\P^*$.    Here we mean {\it triangle\/}
in the sense of the points and lines.
The complement of our line triple in $\P$ is
a union of $4$ open $2$-dimensional
triangular faces.  We ignore these
faces.

Every point of every flat in $X$ has
the structure described above.
In particular, each flat defines a triangle in $\P$
(in the sense above).  Conversely every triangle
in $\P$ comes from a flat in this way.  In short, the flats
are naturally in bijection with the triangles in $\P$.

\newpage

  \section{The Pappus Modular Groups}

\subsection{Actions on the Farey Graph}

Let $\Gamma$ be the set of oriented geodesics in the Farey graph. 
To each directed geodesic in $\Gamma$ we associate the
halfplane it bounds that a person
walking along the geodesic would see on their right.

We now describe $3$ permutations of
$\Gamma$.  Let $e$ be some oriented geodesic in $\Gamma$.
We define $i(e)$ to be the same geodesic but with the opposite
orientation.
The operations $t$ and $b$ are such that
 the geodesics $i(e)$, $t(e)$, and $b(e)$ are the edges of an
  ideal triangle in the Farey graph, and they are oriented
  counterclockwise
  around the triangle.
  Figure 3.1 shows how these operations act on a typical edge $e$.  We also
  shade the associated halfplanes.

 \begin{center}
\resizebox{!}{1.2in}{\includegraphics{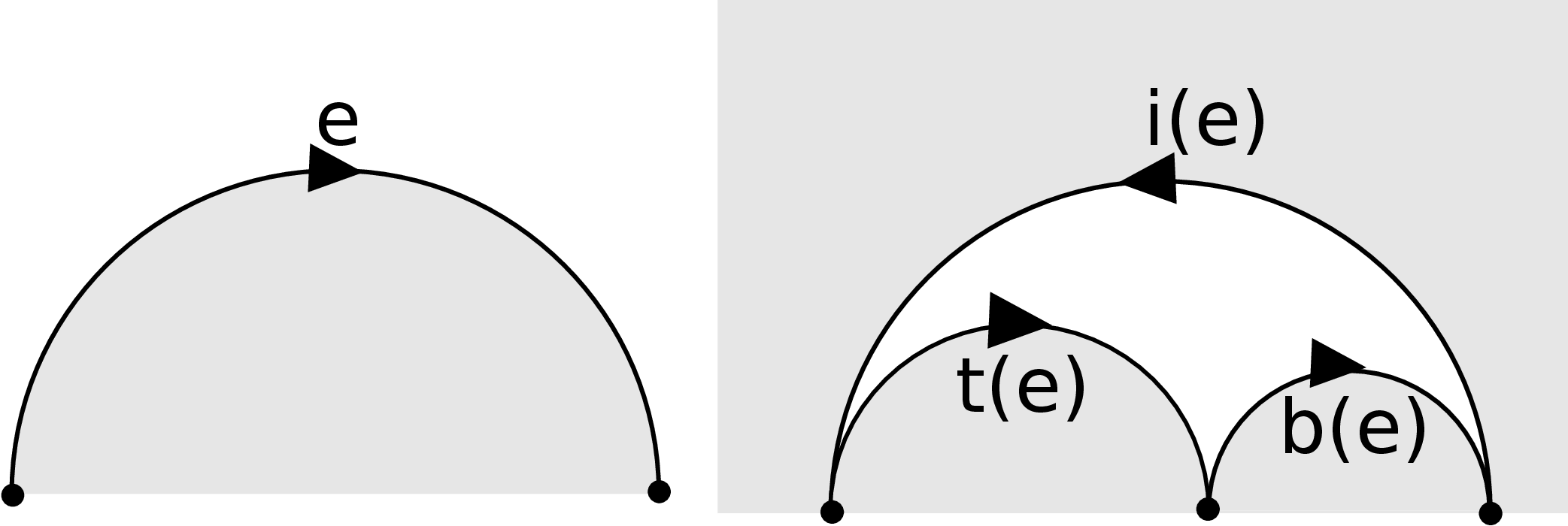}}
\newline
{\bf Figure 3.1:\/} The operations $t$ and $b$.
\end{center}

The halfplanes associated
to $t(e)$ and $b(e)$ are contained in the halfplane associated to $e$.
The halfplanes associated to $i(e)$ and $t(e)$ and $b(e)$ are disjoint.
The operations above satisfy the following
relations.
\begin{equation}
  \label{relations}
  i^2=I; \hskip 10 pt
  tit=b; \hskip 10 pt bib=t, \hskip 10 pt tibi=I, \hskip 10 pt biti=I,
  \hskip 10 pt (it)^3=(ib)^3=I.
\end{equation}
Here $I$ is the identity permutation.

From these relations, and from the nesting properties of the
halfplanes, we see that the group of permutations generated
by $(i,t,b)$ is isomorphic to the modular group, namely the free
product
$\Z/2 * \Z/3$.  Concretely, the group in question is generated by
$i$ and $it$, elements of order $2$ and $3$ respectively.

This action of the modular group on $\Gamma$ is completely
combinatorial.
At the same time, the modular group also acts on $\Gamma$ as
the group of orientation preserving hyperbolic isometries of
$\Gamma$.
This copy of the modular group is generated by an order $3$ rotation
of any
of the ideal triangles, and an order $2$ element which stabilizes 
an edge and reverses its direction.  The combinatorial action
and the geometric action commute with each other.

\subsection{Convex Marked Boxes}

Our description will depart from that in [{\bf S0\/}] though not
in an essential way.  We identify the affine patch of $\P$ with
$\R^2$.  We fix a square $Q \subset \R^2$ whose
sides have unit length and
are parallel to the coordinate axes.  We think of $Q$ as
a solid square.
\newline
\newline
{\bf Model Convex Marked Boxes:\/}
We say that a
{\it model convex marked box\/} is a triple $(Q,t,b)$ where
$t$ is a point in the interior of the top edge of $Q$ and
$b$ is a point in the interior of the bottom edge of $Q$.
We will often abbreviate this terminology to
{\it model box\/}.
\newline
\newline
{\bf Convex Marked Boxes:\/}
We say that a {\it convex marked box\/} is a
triple $(Q',t',b')$ which is projectively equivalent
to some model convex marked box.  That is, there
is a projective transformation
$\Psi$ such that $$(Q',t',b')=\Psi(Q,t,b).$$
We call the edge of $Q'$ containing $t'$ the
{\it top edge\/}. 
We call the edge of $Q'$ containing $b'$ the
{\it  bottom edge\/}.  See Figure 3.2 below.
\newline
\newline
{\bf Reflection Ambiguity:\/}
Note that there are two model convex marked boxes
projectively equivalent to a general
convex marked box.  Let $\rho$ denote the reflection
in the vertical midline of $Q$.  Then $\rho(Q,t,b)$ is another
model convex marked box and $\Psi \circ \rho$ is
another projective transformation, and
$$(Q',t',b')=(\Psi \circ \rho)(\rho(Q,t,b)).$$
We call this the {\it reflection ambiguity\/}.
Up to reflection ambiguity, each convex marked
box is projectively equivalent to a unique
model convex marked box.
\newline
\newline
{\bf The Box Invariant:\/}
We let $x$ denote the distance from $t$ to the
top left vertex of $Q$ and we let $y$ denote the
distance from $b$ to the bottom right vertex of $Q$.
With the reflection ambiguity above in mind, we
write $(x,y) \sim (1-x,1-y)$.
To the marked box $(Q,t,b)$ we associate the
$(\sim)$ equivalence class $[(x,y)]$.
We call this the {\it box invariant\/}.

We extend the box invariant to all other
convex marked boxes by symmetry.  By construction,
two convex marked boxes
are equivalent via a projective transformation
if and only if they have the same
box invariant.  The box invariant can also
be computed using cross ratios, as we discuss
in \S \ref{box}.

\subsection{The Box Operations}

Now we imitate what we did for the oriented Farey graph.

There are $3$ operations on convex marked boxes, which
we call $i$, $t$, and $b$.  Figure 3.2 shows these operations
on a marked box $M$.   The operations $t$ and $b$
encode Pappus's theorem.

 \begin{center}
\resizebox{!}{2.45in}{\includegraphics{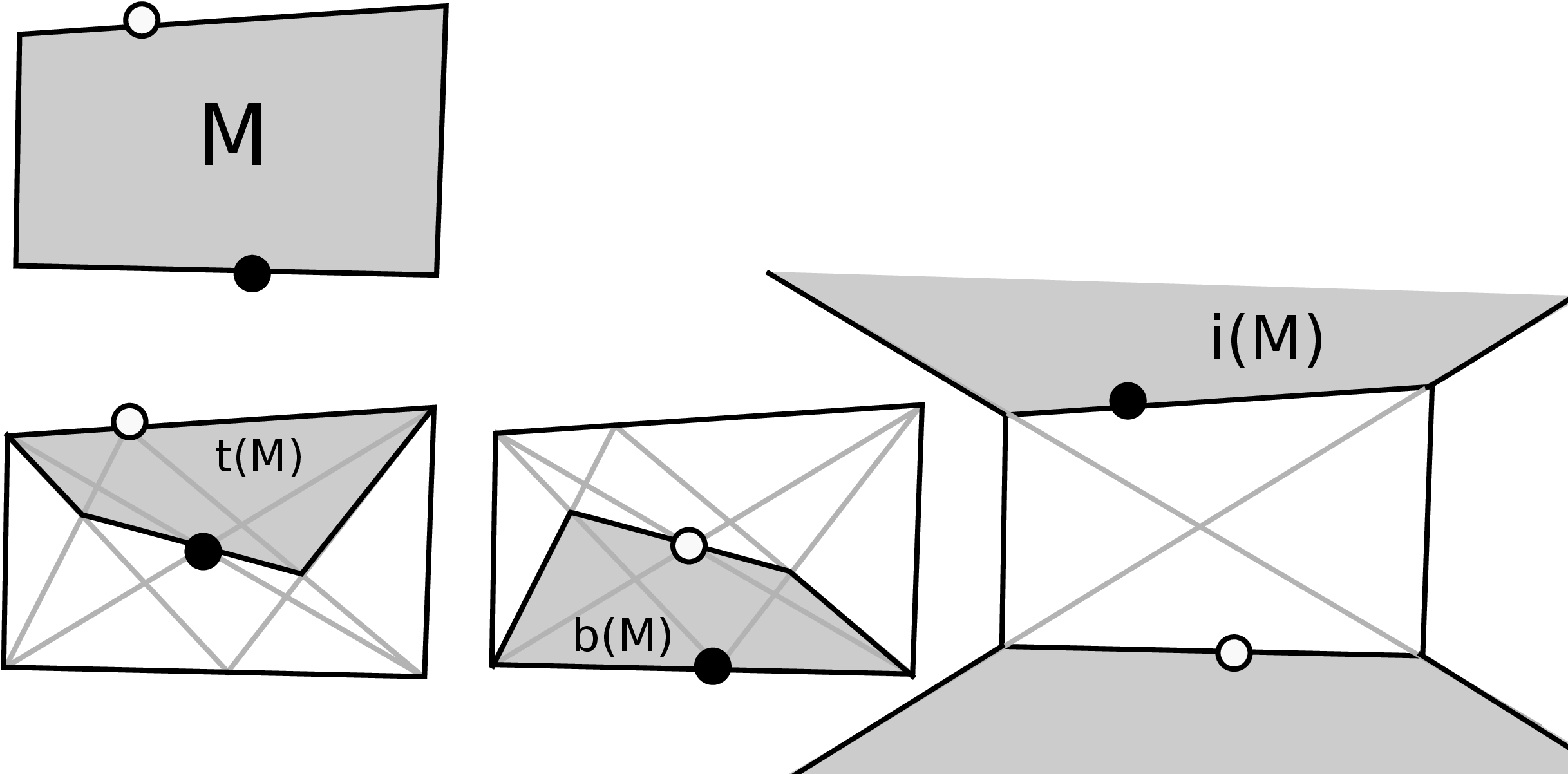}}
\newline
{\bf Figure 3.2:\/} The operations $t,b,i$ on the convex marked box $M$.
\end{center}

On the right hand side, the top of $i(M)$ is the bottom of
$M$ and the bottom of $i(M)$ is the top of $M$.  The top
vertices in each case are colored white and the bottom
ones are colored black.  The box $i(M)$ looks unbounded
in the affine patch, but is indeed a convex marked box
in $\P$.  The shading is meant to go on forever,
so to speak, but we can of course only draw part of it.

The nesting properties of the marked boxes are the same
as for the operations on the oriented Farey graph. Moreover,
one can check easily, by drawing the relevant lines, that
the relations in Equation \ref{relations} hold.  The calculations
are done in [{\bf S0\/}].  From the relations in
Equation \ref{relations} and from the nesting properties,
we see that the group generated by $i,t,b$ is isomorphic
to $\Z/2*\Z/3$ as above.  Thus we have an action of the
modular group on the space of marked boxes.
\newline
\newline
{\bf Remark:\/}
We call this the {\it combinatorial\/} modular group
action, in analogy with the Farey graph case,
but here there is some continuity. The operations
depend continuously (and even algebraically) on the
marked box to which they are applied.

\subsection{A Box Invariant Calculation}
\label{box}

In this section we will start with the model box $M$ with
invariant $[(x,y)]$ and then compute the box invariants for $t(M)$ and $b(M)$.
Given four collinear points
$a,b,c,d \in \P$ we have a {\it cross ratio\/}
\begin{equation}
  [a,b,c,d]=\frac{(a-b)(c-d)}{(a-c)(b-d)}=
  {\rm any\ defined\ entry\ of\/} \frac{(a \times b)(c \times d)}{(a \times c)(b
    \times d)}
\end{equation}
We implement the first formula by identifying the line
containing $a,b,c,d$ with the $X$-axis by
a projective transformation. 
The answer
we get is independent of the choices made.  The second formula
is how we actually do our computations. It
applies when the points $a,b,c,d$ are represented by
vectors in $\R^3$ -- i.e., given in terms of homogeneous
coordinates. We take the $4$ cross products
and then multiply and divide them pointwise.  Each
well-defined entry equals $[a,b,c,d]$.
The cross ratio of $4$ coincident lines has a similar definition.

For the marked box in Figure 3.3 (with corners $s,u,a,c$)
the box invariant is
$[(x,y)]$ where
$$x=[s,t,u,\zeta], \hskip 30 pt
y=[a,b,c,\zeta], \hskip 30 pt \zeta=\overline{su} \cap
\overline{ca}.$$
(In this example $\zeta=[1:0:0]$ but in general it will depend on the
marked box.)
This projectively natural formula works for any marked box.
We name points in Figure 3.3 by vectors in $\R^3$ that represent them.

 \begin{center}
\resizebox{!}{2.2in}{\includegraphics{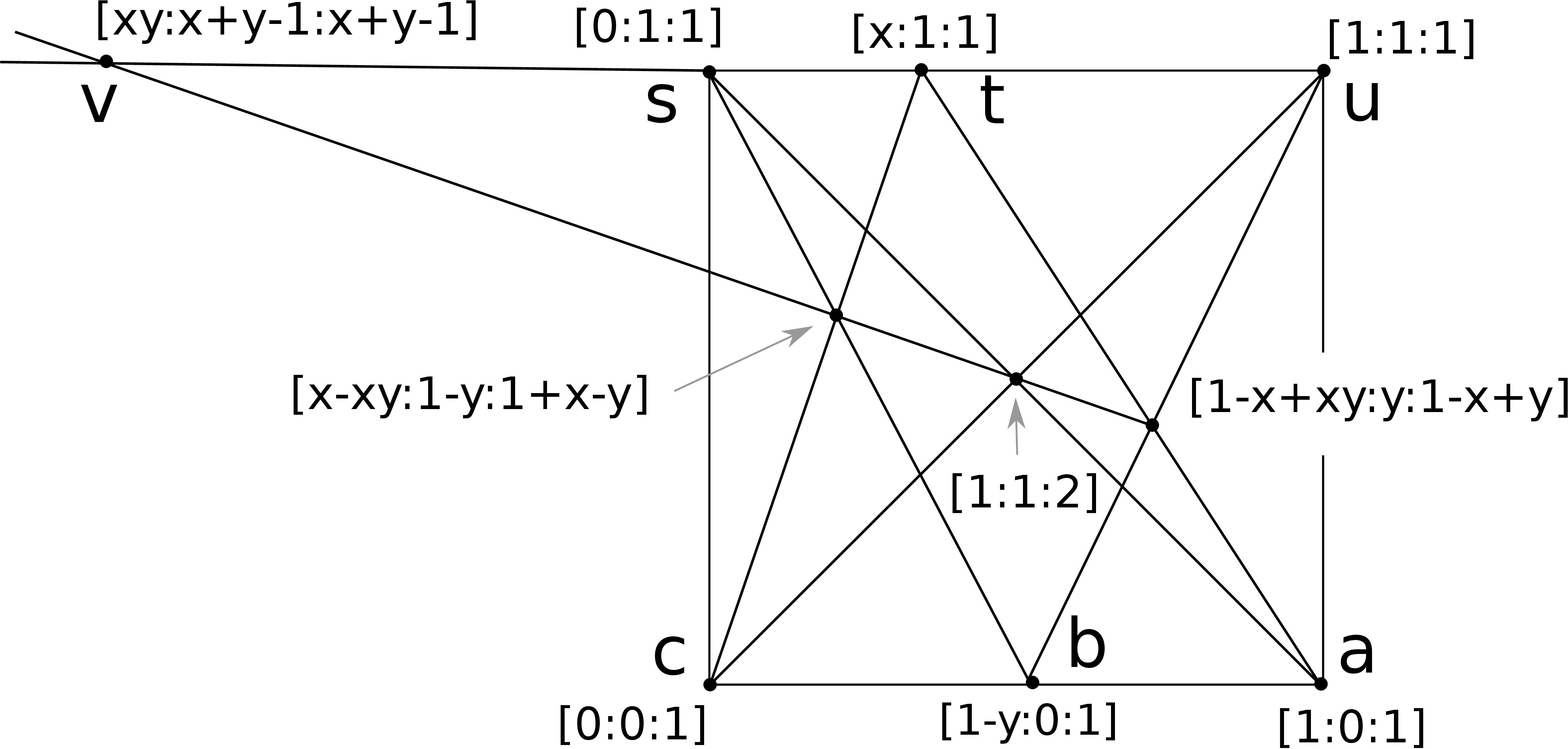}}
\newline
{\bf Figure 3.3:\/}  The relevant coordinates
\end{center}

A calculation shows $t(M)$ and $b(M)$ both have box invariant
$[(1-y,x)]$.  For instance $[s,t,u,v]=1-y$, and
this calculates the first entry in the
box invariant for $t(M)$.

\subsection{Projective Symmetries of an Orbit}

Let $M$ be any marked box.

\begin{lemma}
  \label{order3}
  There is an order $3$ projective transformation $T$ with
the action $t(M) \to b(M) \to i(M)$.
\end{lemma}

\startproof
Since $t(M)$ and $b(M)$ have the same box invariant, there
is a projective transformation $T$ such that
$T(t(M))=b(M)$.
Since the operations commute with projective transformations,
the relations in Equations \ref{relations} give us
$$T(b(M))=T(tit(M))=tiT(t(M))=tib(M)=i(M).$$
$$T(i(M))=T(tib(M))=tiT(b(M))=tii(M)=t(M).$$
The fact that $T^3(M)=M$ forces $T^3$ to be the identity.
\endproof

Lemma \ref{order3} gives us a $\Z/3 * \Z/3$ geometric action, by
projective transformations, on each convex
marked box orbit.  The generators can be taken as $ib$ and $it$.
This group $\langle it,ib\rangle$ has index $2$ in the modular group.
What is missing is a symmetry that maps $M$ to $i(M)$.
Lemma \ref{order3} implies that
$i(M)$ also has box invariant $[(1-y,x)]$.  This is
different from the box invariant $[(x,y)]$ of $M$.
For this reason, $M$ and $i(M)$ are typically
not projectively equivalent.  It turns out that
there is an elliptic polarity
which maps $M$ to $i(M)$.
To make sense of this statement we first need to
interpret marked boxes in a way that is compatible
with the action of dualities.
Our treatment of this topic departs somewhat from that in
[{\bf S0\/}] but not in essential ways.

\subsection{Doppelgangers}

We can specify a convex marked box as a $6$-tuple of
points  $(s,t,u,a,b,c)$.  These points  go in cyclic order
around the boundary of the convex quadrilateral.
Figure 3.4 (left) shows what we mean.
The data $(u,t,s,c,b,a)$ describes the same convex
marked box.  Thus it is really the pair of these
$6$-tuples which describes the marked box.

Given the convex marked box $M$ in $\P$  described by
$(s,t,u,a,b,c)$, we let $M^*$ be the
convex marked box in $\P^*$ described by
$(S,T,U,A,B,C)$.  We call $M$ and $M^*$
{\it doppelgangers\/}.   Figure 3.4 shows
the situation.

 \begin{center}
\resizebox{!}{1.2in}{\includegraphics{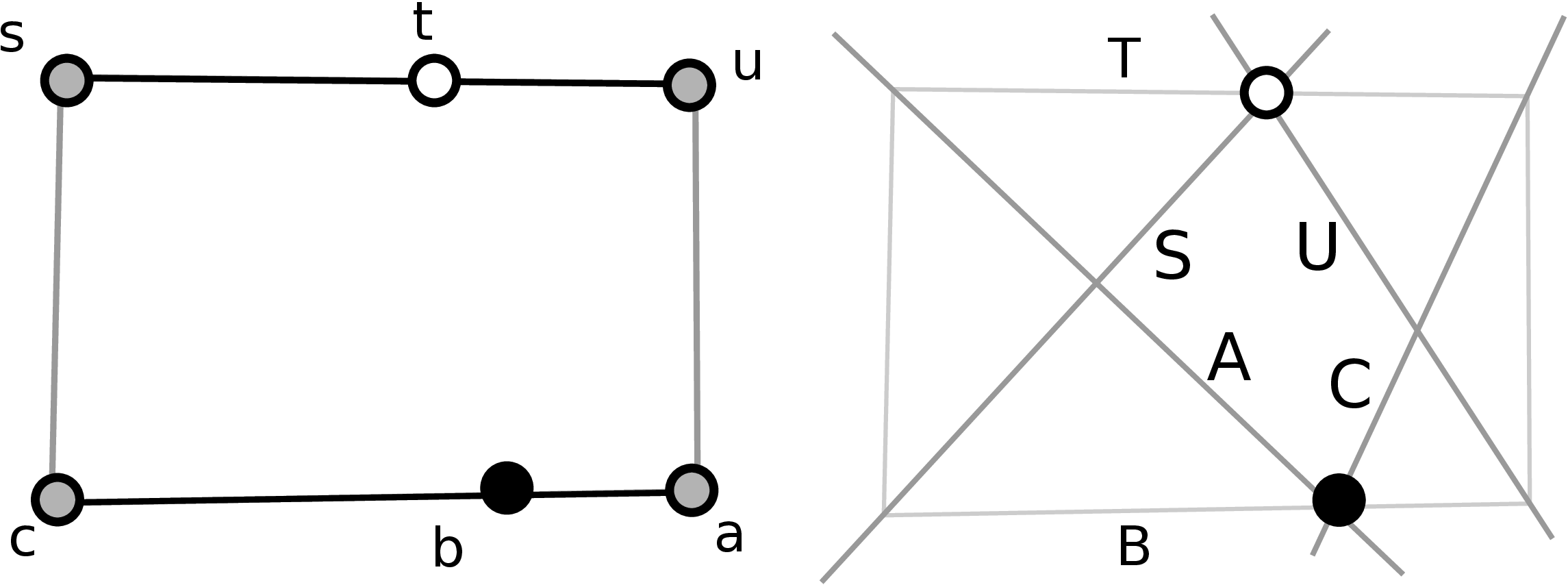}}
\newline
{\bf Figure 3.4:\/} A convex marked box and its doppelganger
\end{center}

It makes sense to perform the operations $i,t,b$ in
$\P^*$.   We just interchange the roles played by
points and lines.   The following result is really
a consequence of the self-dual nature of Pappus's Theorem.

\begin{lemma}[Compatibility]
  Suppose $M$ and $M^*$ are doppelgangers and
  $\gamma$ is one of the three marked box operations.
  Then $\gamma(M)$ and $\gamma(M^*)$ are doppelgangers.
\end{lemma}

\startproof
Consider the operation $i$.
The data for the marked box $i(M)$ is
$(a,b,c,u,t,s)$ and the data for the marked box
$i(M^*)$ is $(A,B,C,U,T,S)$.  In terms of the data,
these are purely combinatorial operations.  From this
we can see that $i(M)$ and $i(M^*)$ are doppelgangers.

Now we consider the operation $t$ (not to be confused with
one of the vertices of our data.)  The data for $t(M)$ is given
by $(s,t,u,a',b',c')$, where
$$a'=\overline{ta} \cap \overline{ub}, \hskip 20 pt
b'=\overline{sa} \cap \overline{uc}, \hskip 20 pt
c'=\overline{tc} \cap \overline{sb}.$$
The data for $t(M^*)$ is given
by$(S,T,U,A',B',C')$, where
$$A'=\overline{(T \cap A)(U \cap B)} \hskip 20 pt
B'=\overline{(S \cap A)(U \cap C)} \hskip 20 pt
C'=\overline{(T \cap C)(S \cap B)}.$$
  Figure 3.5 shows the construction.

 \begin{center}
\resizebox{!}{1.4in}{\includegraphics{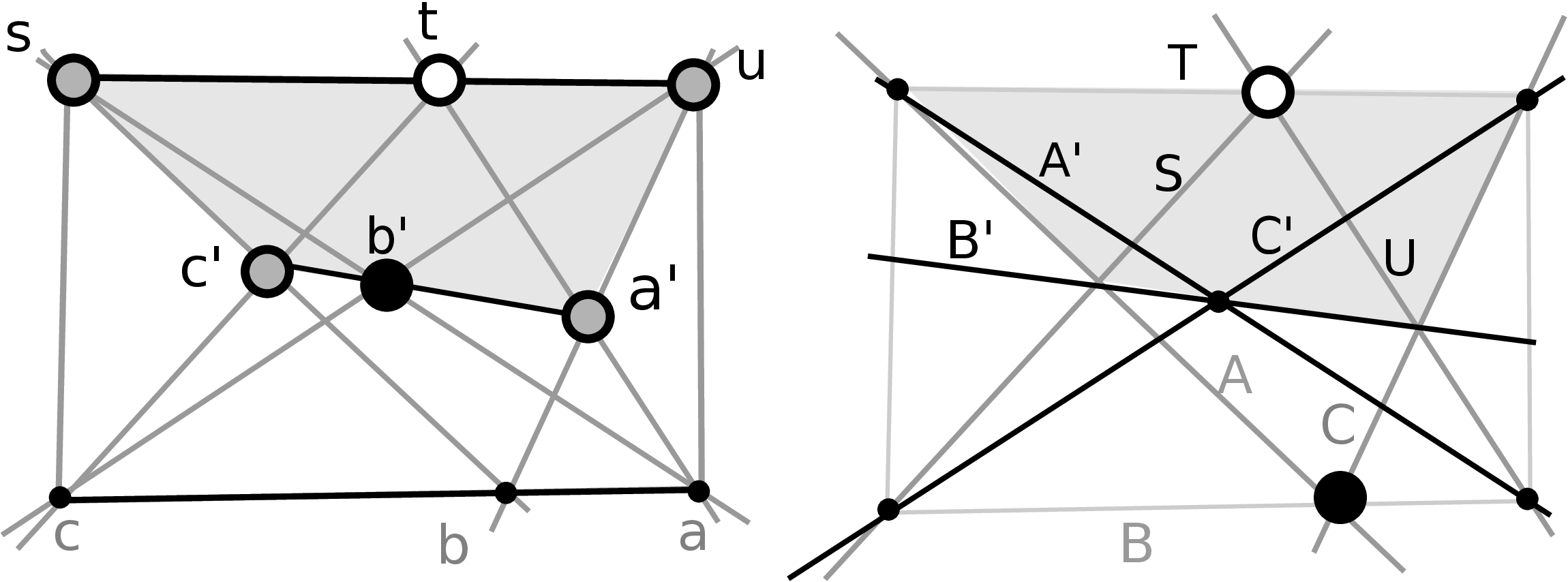}}
\newline
{\bf Figure 3.5:\/} The effect of the operation $t$.
\end{center}

  Looking at these formulas and inspecting the data,
  we see that $t(M)$ and $t(M^*)$ are doppelgangers.
  A similar argument works for the operation $b$.
  \endproof

  \subsection{Duality}

Having defined doppelgangers, we now
discuss how projective dualities interact with
marked boxes.  We keep the notatation from Figure 3.4.

  \begin{lemma}
    \label{pol}
    There is an elliptic polarity which simultaneously
    maps $M$ to $(i(M))^*$ and
    $M^*$ to $i(M)$.
  \end{lemma}

  \startproof
  Notice that $M^*$ and $(i(M))^*$ have the right form
for a marked box in $\P^*$, because they each consist
of $6=3+3$ lines, comprised of two triples
of coincident lines.    Now we make an explicit calculation.
Let $M_1=M$ and $M_2=i(M)$.   We write
$$M_k=(s_k,t_k,u_k,a_k,b_k,c_k)=(u_k,t_k,s_k,c_k,b_k,a_k).$$
The second equality means that the two sets of data are
equally valid representations of $M_k$.
We normalize so that the first sextuple for $M_1$ is represented by vectors
$$
((-1, 1, 0), (x, 1, 0), (1, 1, 0), (1, 0, 1), (y, 0, 1), (-1, 0, 1))
$$
and the sextuple for $M_2$ is represented by vectors
$$
((-1, 0, 1), (y, 0, 1), (1, 0, 1), (-1, 1, 0), (x, 1, 0), (1, 1, 0))
$$
The conditions $|x|,|y|<1$ give convex marked boxes.
The vectors representing the lines of the doppelgangers are as follows.
{\small
$$
S_k=c_k \times t_k, \hskip 5 pt 
T_k=s_k \times u_k, \hskip 5 pt 
U_k=a_k \times t_k, \hskip 5 pt
A_k=s_k \times b_k, \hskip 5 pt 
B_k=a_k \times c_k, \hskip 5 pt 
C_k=u_k \times b_k.
$$
\/}
Here $(\times)$ denotes the cross product.

Let
$$
m=
\left[\matrix{
 1 & -x & -y \cr
 -x & 1 & x y \cr
 -y & x y & 1}\right]
$$
Note that $\det(m)=(x^2-1)(y^2-1)$, so for convex marked
boxes $m$ is nonsingular.
Let $\delta=\Delta \circ m$, where $\Delta$ is the standard polarity.
We compute
$$m(s_1) \times S_2=(0,0,0), \hskip 30 pt (m^{-1})^t(S_1) \times u_2=0.$$
This means that $\delta(s_1)=S_2$ and $\delta(S_1)=u_2$.
Doing $10$ other calculations like this, we see that $\delta$ has the
following action.
$$(s_1,t_1,u_1,a_1,b_1,c_1) \to (S_2,T_2,U_2,A_2,B_2,C_2),$$
$$(S_1,T_1,U_1,A_1,B_1,C_1) \to (u_2,t_2,s_2,c_2,b_2,a_2).$$

Note that $m^t=m$, which means that
$(m^{-1})^t=m^{-1}$.   Recalling that
$\delta=\Delta \circ m$, we compute
$$(\Delta \circ m) \circ (\Delta \circ m) =
\Delta \circ m \circ \Delta \circ m = \Delta \circ \Delta \circ
(m^{-1})^t \circ m = \Delta \circ \Delta \circ m^{-1} \circ m = I.$$
This, $\delta \circ \delta$ is the identity.  Hence
$\delta$ is a polarity.

Finally, when $x=y=0$ we see that $\delta=\Delta$.
In this case $\delta$ is an elliptic polarity.
The set of elliptic polarities is a connected component
of the set of all polarities.  Hence $\delta$ is an
elliptic polarity for all $x,y \in (-1,1)$.
\endproof 

\noindent
{\bf Remarks:\/}
\newline
(1) Our normalization in the Lemma \ref{pol} is
as in [{\bf BLV\/}] rather than as in
[{\bf S0\/}].  Their normalization is very
convenient for the duality calculation.
\newline
(2) Our purely algebraic proof of Lemma \ref{pol}
might leave some readers cold.
In the next chapter we will give a geometric
proof of Lemma \ref{pol}.  The geometric proof
is closer to what we did in [{\bf S0\/}].
\newline
(3)
Our marked box operations  make sense
over essentially any field.  (Fields over $\Z/2$ do not have enough
points on a line to make sense of the constructions.)
In particular, the convexity assumptions
are not really crucial for Equation \ref{relations}.
We needed the convexity to guarantee that
these operations generate a faithful action of the
modular group.  For instance, over finite fields,
the action could not possibly be faithful.
Lemmas \ref{order3} and \ref{pol} also
work over essentially
any field, although in general the notion
of an {\it elliptic\/} polarity does not make sense.
\newline

We say that an {\it enhanced box\/} is a pair
$(M,M^*)$ of doppelgangers. 
If we have a duality $\delta$ and an enhanced
box $(M,M^*)$ then we define
\begin{equation}
  \label{enhanced}
  \delta((M,M^*))=(\delta(M^*),\delta(M)).
\end{equation}
In this way, dualities also act on (enhanced)
marked box orbits.

 Each convex marked box
orbit defines a unique enhanced box orbit on which
the group $\Z/2*\Z/3=\langle i,t,b\rangle$ acts.
If we forget the doppelgangers we just get
the original orbit.  The Duality Lemma
combines with Equation \ref{enhanced}
to give us the missing symmetry.
Now we know that there is a polarity that
interchanges $M$ and $i(M)$ for every
enhanced box $M$ in the orbit.
Combining Lemmas \ref{order3} and \ref{pol}, we get a
geometric modular group action on
each enhanced box orbit.  These are the Pappus
modular groups.

Once we single out some
enhanced box in the orbit we get a natural map
from the oriented geodesics in the Farey graph to the
enhanced box orbit.  The natural map intertwines
the two commuting modular group actions on
the oriented Farey graph and the two commuting
modular group actions on the enhanced box orbit.
The order $3$ isometric rotations of
the ideal triangles in the Farey triangulation correspond
to the projective transformations from Lemma \ref{order3},
and the order $2$ edge stabilizers acting on the
Farey graph correspond to the
polarities from Lemma \ref{pol}.
Theorem \ref{pattern} makes this intertwining map
more canonical and meaningful.

\newpage

  \section{The Symmetric Space Picture}

\subsection{Recognizing the Standard Polarity}

We work in the affine patch, which we identify with $\R^2$.
The standard polarity $\Delta$ maps the point $(r,0)$ to the
line $x=-1/r$.  The reason is that this line consists of
points given in homogeneous coordinates of the form $[-1/r:y:1]$,
and
$$(-1/r,y,1) \cdot (r,0,1)=0.$$
We also note that the action of $\Delta$ on points in the
affine patch is radially symmetric.  For this reason,
Figure 4.1 depicts the action of $\Delta$.
We have $\Delta(p)=L$ and $\Delta(L)=p$.
Here $\|p\|\|q\|=1$.  The action on lines which
do not intersect the unit circle and on points inside
the unit circle satisfies a similar reciprocity rule.
  
 \begin{center}
\resizebox{!}{2.3in}{\includegraphics{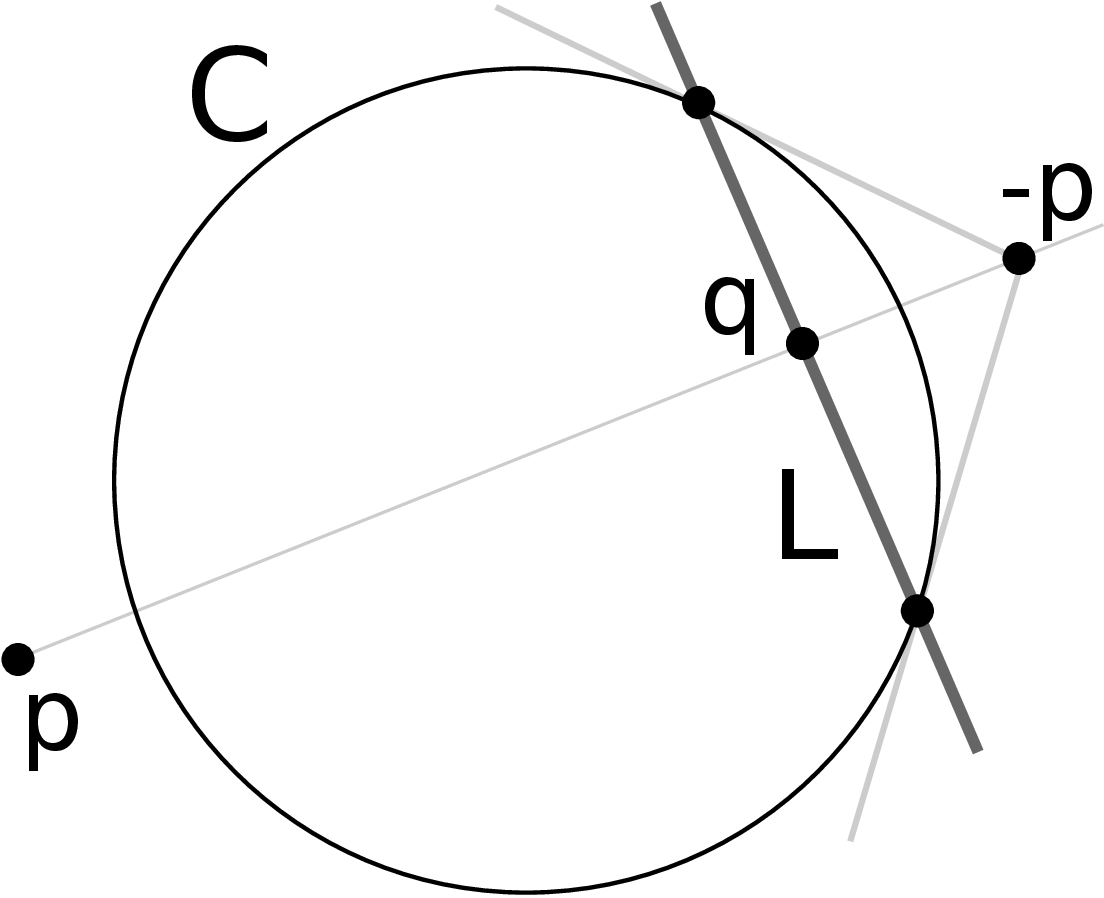}}
\newline
{\bf Figure 4.1:\/} A geometric picture of the standard polarity
\end{center}

There is another way to think about the action of $\Delta$ on
the affine patch that hyperbolic geometers will appreciate.
Think of the open unit disk as the hyperbolic plane in the
Klein model.  In this model, the geodesics are intersections
of straight lines with the open unit disk.
The non-elliptic polarity $\Delta_{\mu}$ represented
by the diagonal matrix $\mu$
with diagonal entries $(-1,-1,1)$ interchanges hyperbolic
geodesics with points according to
the construction shown in Figure 4.1.  The action is
$\Delta_{\mu}(L)=-p$.  At the same time, $\mu$ represents the
the projective transformation $T_{\mu}$ which
is order $2$ rotation about the origin in $\R^2$.
Figure 4.1 depicts the composition $T_{\mu} \circ \Delta_{\mu}$,
and this polarity is represented by $\mu^2$, the identity matrix.
Hence  $\Delta=T_{\mu} \circ \Delta_{\mu}$.

\subsection{Lemma \ref{pol} revisited}

Here we give a second proof of Lemma \ref{pol}.
Figure 4.2, an elaboration of
[{\bf S0\/}, Figure 2.4.2],
shows how we normalize the marked box $M$ so that
$\Delta$ maps $M$ to $i(M)$.  The lines $T$ and $B$ go through
the origin in $\R^2$.  The points $t$ and $b$ lie at $\infty$ on these
lines.  We also have $\|a\|\|c\|=\|u\|\|s\|=1$.
The shaded regions show the convex quadrilaterals
(in $\P$) comprising $M$ and $i(M)$ respectively.

 \begin{center}
\resizebox{!}{4in}{\includegraphics{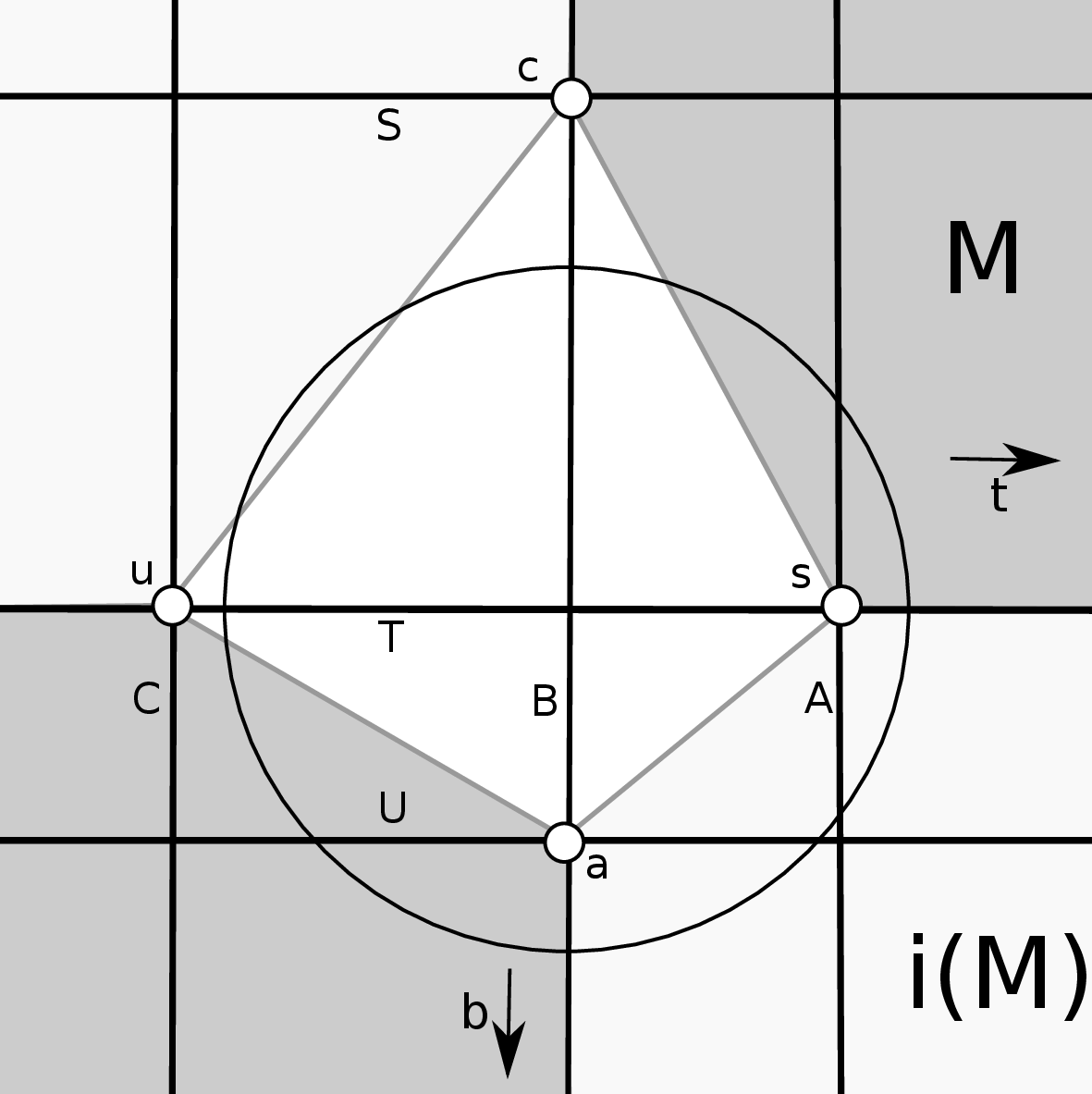}}
\newline
{\bf Figure 4.2:\/} A normalized view of $M$ and $i(M)$.
\end{center}

Comparing Figures 4.1 and 4.2 we see that $\Delta$ has the action
$$(s,t,u,a,b,c) \to (C,B,A,S,T,U)=i(S,T,U,A,B,C).$$
Thus $\Delta$ maps $M$ to $i(M)$ as claimed.

By symmetry, there is an elliptic polarity which maps each
marked box $M$ to the marked box $i(M)$.

\subsection{The Farey Patterns}

Let $\cal M$ be a marked box orbit and let
$\Lambda_{\cal M}$ be the Pappus modular group of projective
symmetries of $\cal M$.  For each marked box
$M$ of $\cal M$ we have an elliptic polarity
$\delta_M \in \Lambda_{\cal M}$ which swaps $M$ and $i(M)$.
We let $p_M \in X$ be the fixed point of
$\delta_M$.   Equivalently, $p_M$ is the
unit ellipsoid for the quadratic form which
gives $\delta_M$.

Note that $p_{i(M)}=p_M$.
The points
in $X$ we have just constructed correspond
to the fixed points of the order $2$ elements
of $\Lambda_{\cal M}$.
There is also a flat $f_M$ we can associate to the
marked box $M$.   Let $(t,T)$ be the top flag
of $M$ and let $(b,B)$ be the bottom flag.
We have a triangle in $\P$ whose vertices are
$(t,b,T \cap B)$ and whose sides are
$(T,B,\overline{tb})$.   The flat $f_M$ is the
flat associated to this triangle, as discussed in
\S \ref{symm}.

\begin{lemma}
We have $p_M \in f_M$.
\end{lemma}

\startproof
By symmetry, it suffices to consider the
case when $M$ is normalized as in Figure 4.2.
In this case $\delta_M=\Delta$ and
$p_M$ is the origin in $X$, the point which
names the unit ball in $\R^3$.
The triangle is the one defined by the coordinate
axes and coordinate planes. Thus the flat $f_M$
is the standard flat.  In this case
$f_M$, the standard flat,  contains
the origin $p_M$ of $X$.
\endproof

We can also associate to $M$ a $1$-parameter subgroup of
$SL_3(\R)$.   Consider the set of matrices whose eigensystem is
\begin{equation}
  \label{diag}
(\lambda,1/\lambda,1), \hskip 30 pt (t,b,T \cap B), \hskip 30 pt
\lambda \in \R.
\end{equation}
This set forms a $1$-parameter subgroup $H_M$ of $SL_3(\R)$
which acts as translation on $f_M$.  When
$M$ is normalized as in Figure 4.2, the matrices are
diagonal, with diagonal entries $\lambda,1/\lambda,1$.
The orbits of $H_M$ in $f_M$ give a foliation of $f_M$ by
medial geodesics which limit on one end to the flag
$(t,T)$ and on the other end to $(b,B)$.  There is a unique
member of this foliation which contains $p_M$.  This is
the medial geodesic we associate to $M$.  We call it
$\gamma_M$.  We choose the orientation so that
$\gamma_M$ is asymptotic to $(t,T)$ in the backwards
direction and asymptotic to $(b,B)$ in the forwards direction.
\newline
\newline
{\bf Definition:\/}
Let $\Gamma_{\cal M\/}$ denote the union of all the oriented
geodesics of the form $\gamma_M$ for $M \in \cal M$.
\newline

This is our pattern of geodesics in $X$ corresponding
to the Pappus modular group $\Lambda_{\cal M}$.
Since the geodesics in $\Gamma_{\cal M}$ are bijectively
associated with marked boxes in $\cal M$, it makes
sense to speak of the combinatorial modular group
action on $\Gamma_{\cal M}$.

Once we have chosen a single oriented geodesic
$\gamma \in \Gamma$ and a single oriented
geodesic $\gamma_M \in \Gamma_{\cal M\/}$ we have a
canonically defined bijection $\Gamma \to \Gamma_{\cal M}$ which
has the following virtues:
\begin{itemize}
\item The bijection is an isometry on each
  geodesic.
\item The bijection intertwines the two combinatorial
  modular group actions.
\item The bijection intertwines the two geometric
  modular group actions.
\end{itemize}

\begin{lemma}
  Two geodesics in $\Gamma_{\cal M}$ are asymptotic
  if and only if the corresponding geodesics are asymptotic
  in $\Gamma$.
\end{lemma}

\startproof
Given the nesting properties of the convex quads underlying
our marked boxes in $\cal M$ we see that two flags
associated to boxes in $\cal M$ coincide if and only if
the corresponding rational points in $\R \cup \infty$
are the same.  Our lemma follows from this fact, and
from the construction of the geodesics in $\Gamma_{\cal M}$.
\endproof

One of the basic results we proved in
[{\bf S0\/}] is that the flags associated to the marked boxes in
$\cal M$ are dense in a continuous loop in $\cal P$.
Typically this is a fractal loop, as Figure 1.2 suggests.
Our bijection between geodesics in the Farey triangulations
and geodesics in $\Gamma_{\cal M}$ induces a homeomorphism
from $\R \cup \infty$ to this loop.

There is one
  special case where
  the underlying marked boxes are {\it symmetric\/}.
  In this case the marked box has fixed by an order $2$
  reflection which fixes a point $p$ and a line $\ell$.
  The box invariants in this case are
    $(1/2,1/2)$.  In this case, the Pappus
  modular group stabilizes a totally geodesic
  copy of $\H^2$ sitting in $X$.  The
  corresponding loop in $\cal P$ consists of all the
  flags whose points lie on $\ell$ and whose lines contain $p$.
  
  In the symmetric case,
  the  corresponding Farey pattern is isometric to the
  usual Farey triangulation.    This symmetric
  case just recreates an isometric copy of the
  Farey pattern in $X$.  As we move away from
  the symmetric case, we get nontrivial
  deformations of the Farey pattern.

\subsection{Symmetry}

In this section we analyze the symmetry of
our Farey pattern.  Let $\Lambda_{\cal M}$ denote
the Pappus modular group associated to the
marked box orbit $\cal M$.
Let $\widehat \Lambda_{\cal M}$ denote the
symmetry group of the Farey pattern
$\Gamma_{\cal M}$.   This is the subgroup
of isometries which preserve the pattern.

\begin{lemma}
  $\Lambda_{\cal M} \subset \widehat \Lambda_{\cal M}$.
    \end{lemma}

  \startproof
Note that $\gamma_M$ and $\gamma_{i(M)}$ are
the same geodesic but with opposite orientations.
$\delta_M$ is an isometry of the underlying geodesic
and just reverses the orientation. Thus
$\delta_M$ swaps $\gamma_M$ and $\gamma_{i(M)}$.
Given the naturality of our construction, $\delta_M \in \widehat
\Lambda_{\cal M}$.
Also given the naturality of our construction, the
order $3$ projective transformation $T$ whose
orbit is $i(M) \to t(M) \to b(M)$ preserves
the corresponding union of $3$ oriented geodesics
in $\Gamma_{\cal M\/}$.
The same naturality shows that $T \in \widehat \Lambda_{\cal M}$.
Since $\Lambda_{\cal M}$ is generated by the
order $2$ and order $3$ elements, $\Lambda_{\cal M} \subset
\widehat \Lambda_{\cal M}$.
\endproof

\begin{lemma}
  For a generic marked box orbit $\cal M$, we have
  $\Lambda_{\cal M}=\widehat \Lambda_{\cal M}$.
\end{lemma}

\startproof
A {\it Farey triangle\/} in $\Gamma_{\cal M}$ is a
triple of unoriented medial geodesics in the Farey pattern
which correspond to the boundary of an
ideal triangle in the Farey triangulation.
Given the asymptotic
properties of geodesics in the Farey pattern,
every isometry of $\Gamma_{\cal M}$ permutes
the Farey  triangles.

Suppose we have an isometry $\psi$ of $X$ which preserves
$\Gamma_{\cal M}$.  Since $\Lambda_{\cal M}$ acts transitively
on the oriented medial geodesics in our pattern, we can assume
without loss of generality that $\psi$ stabilizes some
oriented medial geodesic $\gamma_{i(M)}$ in the pattern.
Since $\Lambda_{\cal M}$ has a polarity reversing the
orientation of this geodesic, we can assume without
loss of generality that $\psi$ is a linear transformation
which stabilizes the unoriented geodesic in the pattern
which is asymptotic to the two flags $(t,T)$ and $(b,B)$.
Here $(t,T)$ and $(b,B)$ respectively are the top and bottom
flags associated to the marked box $M$.

Consider the fourth power $\psi^4$.  This element must
stabilize the oriented geodesic $\gamma_M$ and also
each of the two Farey triangles which involve
$\gamma_M$.   In terms of the projective action of
$\psi^4$ on $\P$, this element must
fix the $3$ points $(t,b,T \cap B)$ and also the
$4$th point $w_1$ which is common to
both $b(M)$ and $t(M)$.  The point I mean is
simultaneously the top point of $b(M)$ and the bottom point of $t(M)$.
Generically these $4$ fixed points are in general position,
and this forces $\psi^4$ to be the identity.

To get a finer analysis, we
normalize $M$ as in Figure 4.2.
Once we do this, we see that $\psi$ permutes two points
at infinity (meaning in $\P-\R^2$) and hence $\psi$ preserves
the affine patch $\R^2$.  Also, $\psi$ permutes the coordinate
axes in $\R^2$.  These facts, together with the fact that
$\psi$ has finite order, forces $\psi$ to act as an isometry of
$\R^2$.

 \begin{center}
\resizebox{!}{3.4in}{\includegraphics{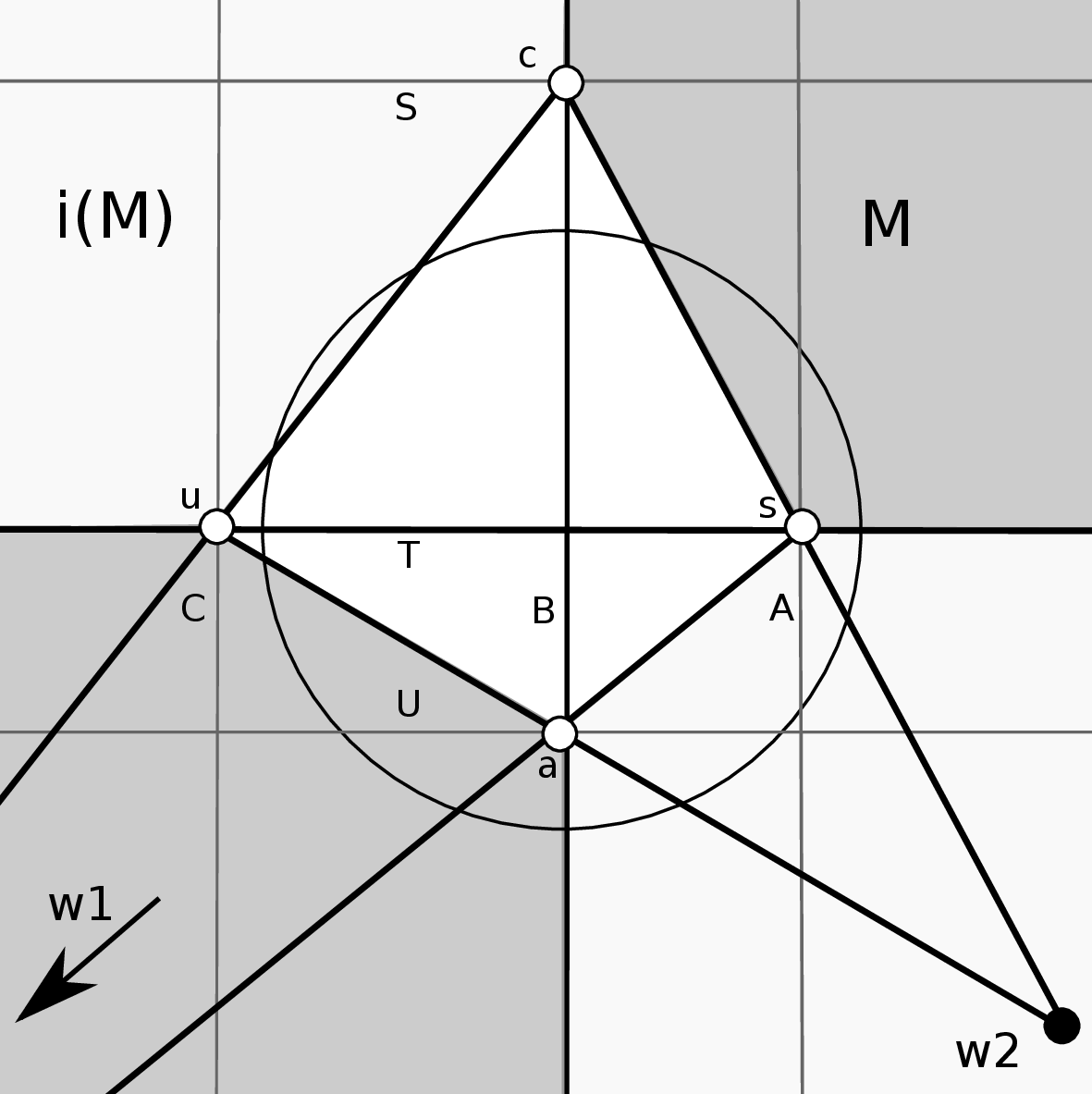}}
\newline
{\bf Figure 4.3\/} The points $w_1$ (off screen) and $w_2$.
\end{center}

Consider the two points
$$w_1=\overline{cu} \cap \overline{as}, \hskip 30 pt
w_2=\overline{au} \cap \overline{cs}.$$
The point $w_1$ is the point common 
to $b(M)$ and $t(M)$ and the point $w_2$ is the
point common to $bi(M)$ and $ti(M)$.  These
points correspond to the limit points of the two
abovementioned Farey triangles which are
opposite $\gamma_M$.   The element $\psi$
must either preserve or permute these two
points.  Generically, they are not the same distance
from the $X$-axis (or from the $Y$-axis) and
so $\psi$ must preserve them.  But then the
same argument as above shows that $\psi$ is
the identity.  So, for a generic choice of $\cal M$
we have $\Lambda_{\cal M}=\widehat \Lambda_{\cal M}$.
\endproof

\subsection{Disjointness}

To finish the proof of
Theorem \ref{pattern} we need to
show that every two distinct geodesics
$\gamma_{M_1}$ and $\gamma_{M_2}$ in
our pattern $\Gamma_{\cal M}$ are
disjoint.   We write $\gamma_1=\gamma_{M_1}$, etc.
We will prove the stronger result that 
the flats $f_1$ and $f_2$ are disjoint in all cases.

\begin{lemma}
  \label{special}
  If $M_1$ and $M_2$ share a common flag then
  the flats $f_1$ and $f_2$ are disjoint.
\end{lemma}

\startproof
Suppose that $f_1$ and $f_2$ intersect.
We can normalize so that this intersection point is
the origin in $X$.   Then $f_1$ and
$f_2$ are two geodesics through the origin
that define the same point on $\partial X$ in one
direction. This forces $f_1 \cap f_2$ to contain a
medial geodesic.   But then both flats are determined
by the same pair of flags and hence coincide.
We are really invoking the more general principle
that distinct flats can only share a geodesic if
it is one of the singular geodesics.
\endproof

It remains to consider the case
when $M_1$ and $M_2$ have no flags in common.
The convex quadrilaterals underlying
$M_1$ and $M_2$ are, like the halfplanes
associated to the Farey triangulation, either
nested or disjoint.  Since we have
$f_{M_j}=f_{i(M_j)}$ we can replace
$M_j$ by $i(M_j)$ if necessary and then
assume that
$M_1$ and $M_2$ have nested underlying
convex quadrilaterals $Q_1$ and $Q_2$.

Suppose for the sake of contradiction that $f_1$ and $f_2$
intersect.  We can normalize so that the intersection point is
the origin in $X$ and that $f_1$ is the standard flag.  We
can further normalize so that
$H_1$ is the subgroup of diagonal matrices with
entries $(\lambda,1/\lambda,1)$.
 The box $M_1$ almost looks like the box
$M$ in Figures 4.2 and 4.3.  The only difference is that the
reciprocity condition is gone.  That is, $M_1$ is the image
of $M$ under a positive diagonal linear transformation acting on $\R^2$.
The important thing is that the top and bottom lines of
$M_1$ are the coordinate axes in $\R^2$.

Consider the triangle $\tau_2$ associated to $f_2$.
Since $f_2$ contains the origin in $X$, the 
triangle $\tau_2$ is {\it orthogonal\/} in the sense
that the three $1$-dimensional subspaces
representing its vertices are mutually orthogonal.
Equivalently, $\tau_2$ is invariant under the
standard polarity $\Delta$. This is to say that
$\Delta$ maps each point of the triangle to the
opposite line defined by the triangle.

Define the {\it positive cone\/} in $\R^2$ to be
the interior of the union of the $++$ and $--$ quadrants.
The important thing about $M_1$ is that its top and bottom
lines cut out the positive cone.

\begin{lemma}
  \label{triangle}
  Two of the lines defined by $\tau_2$ have negative slope.
  \end{lemma}

  \startproof 
  Any convex quadrilateral
strictly contained in the convex quad $Q_1$ associated to $M_1$ lies
in the positive cone.  Since $M_1$ and $M_2$ have
no flags in common, and $Q_2 \subset Q_1$, we see that
$Q_2$ lies in the positive cone.
This means that the top and
bottom points of $Q_2$ are contined in the positive cone.
  
  Let $p$ be one of the points of $\tau_2$ in the positive
  quadrant.  The line of $\tau_2$ opposite to $p$ is the
  image $\Delta(p)$.  Given the action of $\Delta$
  depicted in Figure 4.1, we see that $\Delta(p)$ has
  negative slope.   The extreme cases happen when
  $p$ is near one of the two coordinate axes.  In these
  cases, the slope of $\Delta(p)$ is near $0$ or $-\infty$.
  So, if $p_1$ and $p_2$ both lie in the positive quadrant
  then $\Delta(p_1)$ and $\Delta(p_2)$ have negative slope.
  \endproof

  The top line and bottom line of $M_2$ (meaning the lines that
  extend the top and bottom edges of $Q_2$) are two of
  the three lines of $\tau_2$.
  By Lemma \ref{triangle}, at least one of them has negative slope.
  We will contradict this by showing that both these lines have
  non-negative slope.
    
  Let $S_1$ denote the union of the two edges of the
  $Q_1$ which are not the top or bottom edge.  The
  sides of $S_1$ are the side edges.  

  \begin{lemma}
    The top and bottom lines of $M_2$ intersect both segments of $S_1$.
  \end{lemma}

  \startproof
  We say that a line $\ell$ is {\it adapted\/} to a marked box $M$
  if it intersects both the side edges of the associated convex quad $Q$.
  Looking at Figure 3.2, we see that any line adapted to $t(M)$
  is also adapted to $M$.  Likewise, any line adapted to $b(M)$ is
  also adapted to $M$.  Iterating, we see that any line adapted
  to any marked box in the semigroup orbit
  ${\cal O\/}=\langle r,b \rangle (M)$
  is adapted to $M$.  But any marked box in $\cal M$ whose associated
  convex quadrilateral is contained in $Q$ lies in $\cal O$.
  Applying this argument to $M=M_1$ we see that every line
  adapted to $M_2$ intersects both sides of $S_1$.  
  This result applies to the top and bottom lines of $M_2$, which
  are adapted to $M_2$.
  \endproof

  When $M_1$ is exactly as in Figure 4.2, the two segments of $S_1$,
  namely $\overline{cs}$ and $\overline{au}$ both have their
  endpoints on opposite boundary components of the positive cone.
  The same thing therefore holds for an image of $M$ under a
  positive diagonal matrix.  But then any line which hits both
  edges of $S_1$ has non-negative slope.  Hence the top line and
  bottom line of $M_2$ both have positive slope.
  This is a contradiction.
  Our proof of Theorem \ref{pattern} is done.

\newpage
  \section{Prisms and Bending}

\subsection{Triple Products of Flags}

Suppose that $J= \{(p_k,\ell_k)\}$ is a triple of general position flags.
We represent $(p_k,\ell_k)$ by a pair of vectors $(P_k,L_k)$ where
$P_k \cdot L_k=0$.   The {\it triple product\/} is defined as
\begin{equation}
  \label{tripp}
 \xi(J)= \frac{(P_1 \cdot L_2)(P_2 \cdot L_3)(P_3 \cdot L_1)}{(P_2 \cdot L_1)(P_3 \cdot L_2)(P_1\cdot L_3)}.
\end{equation}
This is an invariant that is independent of the way we represent our
flags with vectors.
Compare [{\bf FG\/}].
If $J'$ is a permutation of $J$ then $\chi(J')=\chi(J)$ or $\chi(J')=1/\chi(J)$, depending
on whether we use an even or odd permutation.

When the triple product is not equal to $-1$, two
triples of flags are equivalent under a projective transformation
if and only if they have the same triple product.
Indeed, when the triple of lines of the flag are not coincident -- and
this is generically the case -- we can normalize a triple of flags as in Figure 5.1.
This picture has $3$-fold rotational symmetry.  The invariant tells
``how far out'' the points are relative to the central triangle.
Compare [{\bf S0\/}, Figure 2.4.1]

 \begin{center}
\resizebox{!}{2.3in}{\includegraphics{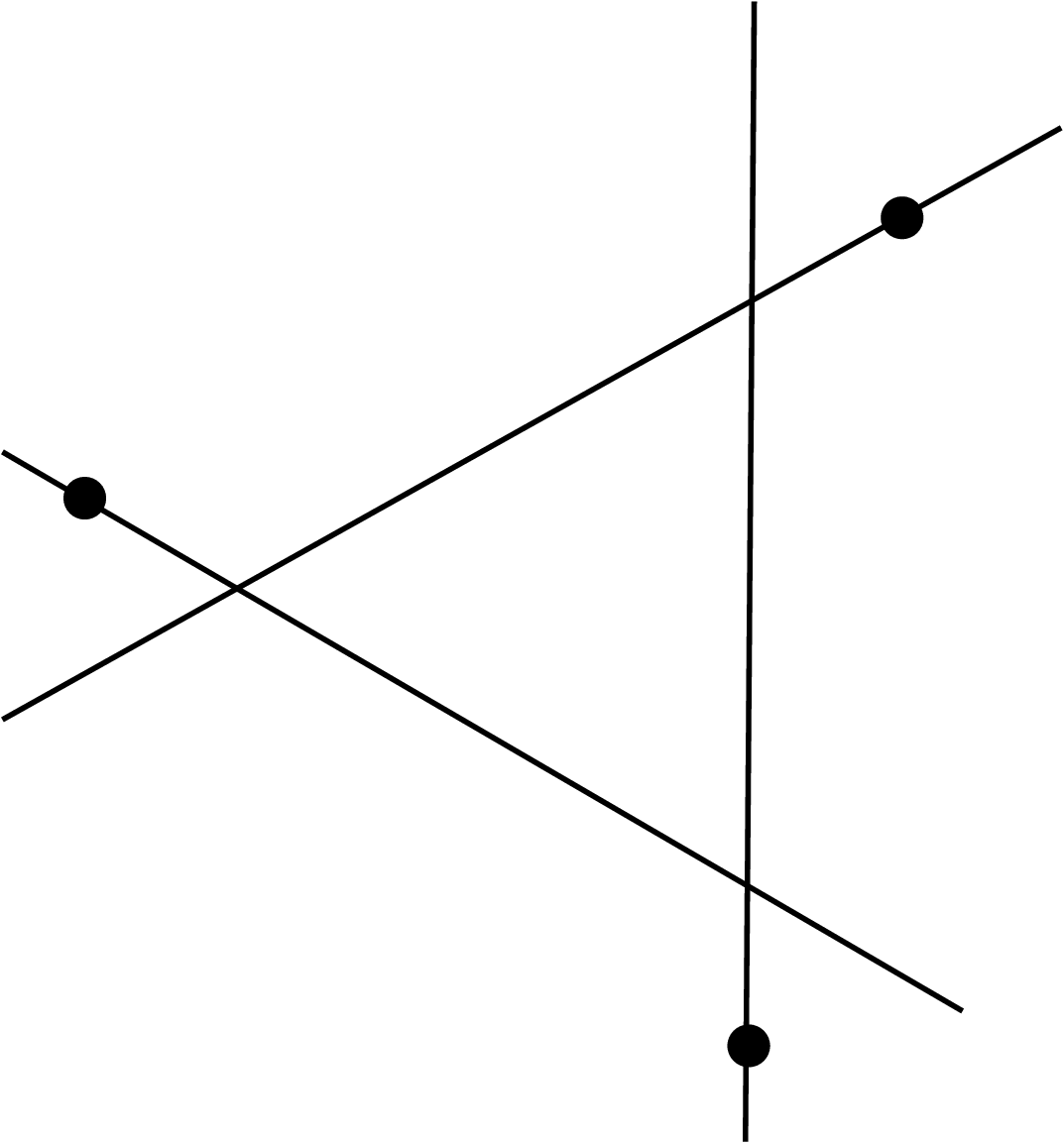}}
\newline
{\bf Figure 5.1\/} A normalized triple of flags
\end{center}
   
For the picture shown the triple product is negative.  The positive case
happens when the points are in the interiors of the edges of the
central triangle.
From this normalized structure we see that
the only projective transformations stabilizing
the generic flag is an order $3$ rotation. Compare
Lemma \ref{order3}.  Lemma \ref{triple} below
discusses how polarities interact with triples of flags.

\subsection{The Triple Invariant and the Character Variety}

We define the triple product of a marked box $M$ to be the triple product of the
flags $$\tau(i(M)),\  \tau(t(M)),\   \tau(b(M)).$$  Here
$i,t,b$ are our three box operations and $\tau$ means ``take the top flag
associated to the marked box''.  Our definition favors the top over the bottom.
Were we to use the bottom flags, we would get the reciprocal.

Now we compute the triple product in terms of the box invariants.
Suppose our marked box has box invariant $[(x,y)]$. Using the
coordinates given in Figure 3.3, we compute
\begin{equation}
  \chi(M)=-\frac{x(1-x)}{y(1-y)} <0
\end{equation}
This formula respects our equivalence relation
$(x,y) \sim (1-x,1-y).$
Since the box invariant for $i(M)$ is
$[(1-y,x)]$ we see that
$\chi(i(M))=1/\chi(M)$.  This means that,
up to taking reciprocals, $\chi(M)$ is independent
of which $M$ we take within a marked box orbit
$\cal M$.

We define
the {\it triple product invariant\/} of the Pappus modular group
$\Lambda_{\cal M}$ to be
\begin{equation}
  \chi(\Lambda_{\cal M}) = |\log(-\chi(M))|.
\end{equation}
This is well-defined independent
of choices because $|\log(r)|=|\log(1/r)|$.
Here is a contour plot of the triple invariant as a
function of the box invariant.

 \begin{center}
\resizebox{!}{1.6in}{\includegraphics{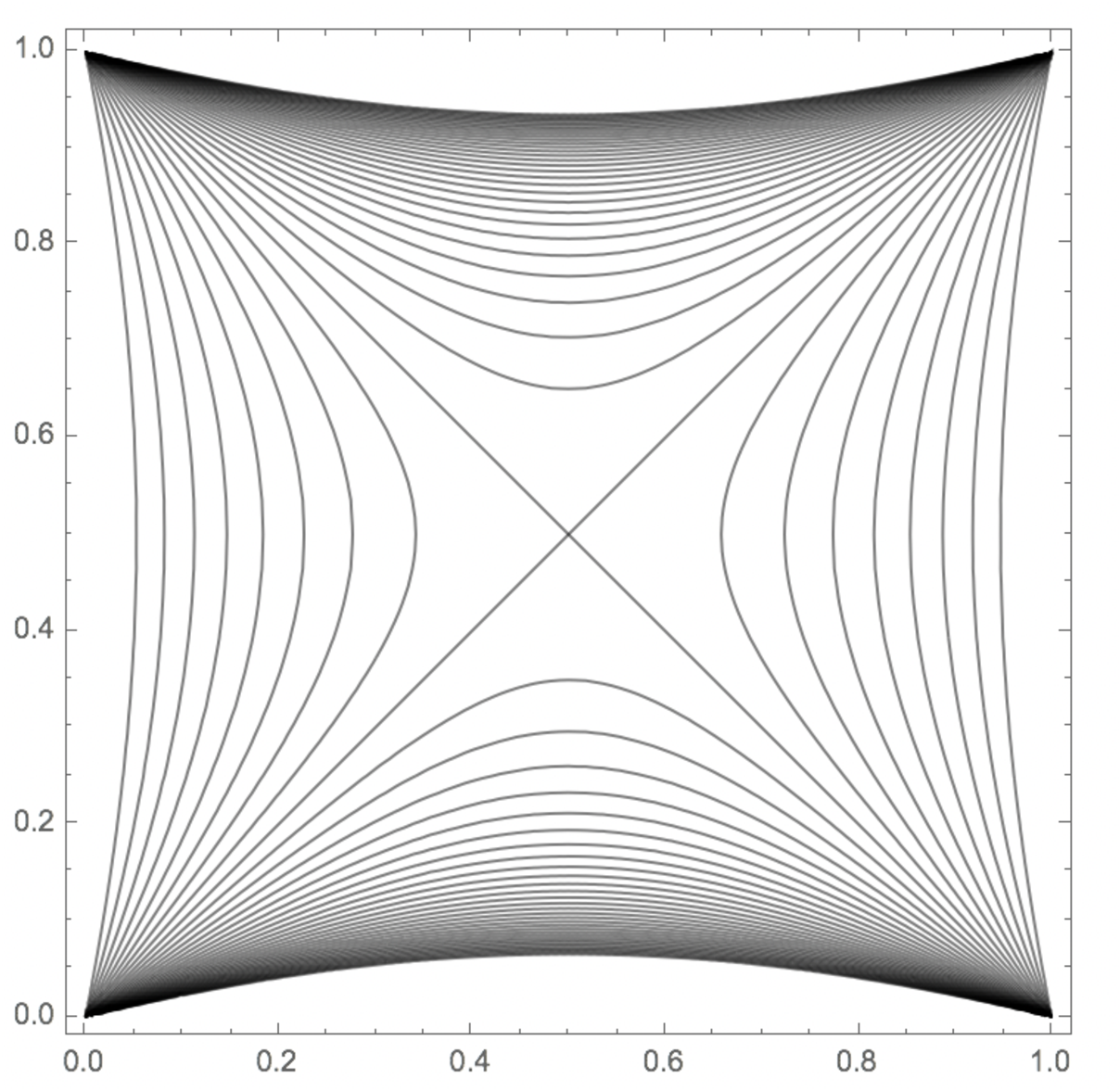}}
\newline
{\bf Figure 5.2\/} Level sets of the triple invariant
\end{center}

The {\it character variety\/} is the space of
Pappus modular groups modulo conjugation.
The conjugation can be by isometries coming
from projective transformation or by those
coming from dualities.

To specify an element of the character variety,
we choose a marked box orbit.  There are
$2$ marked box invariants associated to this
orbit, namely
$\{\rho^k((x,y))|\ k=0,1,2,3\}$ where
$\rho$ is the order $4$ rotation
fixing $(1/2,1/2)$.

It is also worth mentioning that
there is a projective
transformation which maps the box
with invariant $[(x,y)]$ to the box with
invariant $[(y,x)]$ but switches the top
edge and the bottom edge.  This projective
transformation does not completely
respect the markings of the box but
the only thing it does is switch the top
and the bottom.  The two Pappus modular
groups are the same after we change the
names of the generators.  We will consider
the two {\it representations\/} distinct.

This is all the symmetries we have.  So, the
character variety is an open cone: the quotient
$C=(0,1)^2/\langle \rho \rangle$.  Referring to
Figure 5.2, one can take a fundamental domain
for the cone to be the upper quadrant shown
in Figure 5.3.   One gets the cone by
identifying the bounding diagonals in Figure 5.3.
Reflection in the vertical midline
is an involution which, as mentioned above,
preserves the Pappus modular group but
changes the representation.

 \begin{center}
\resizebox{!}{1.2in}{\includegraphics{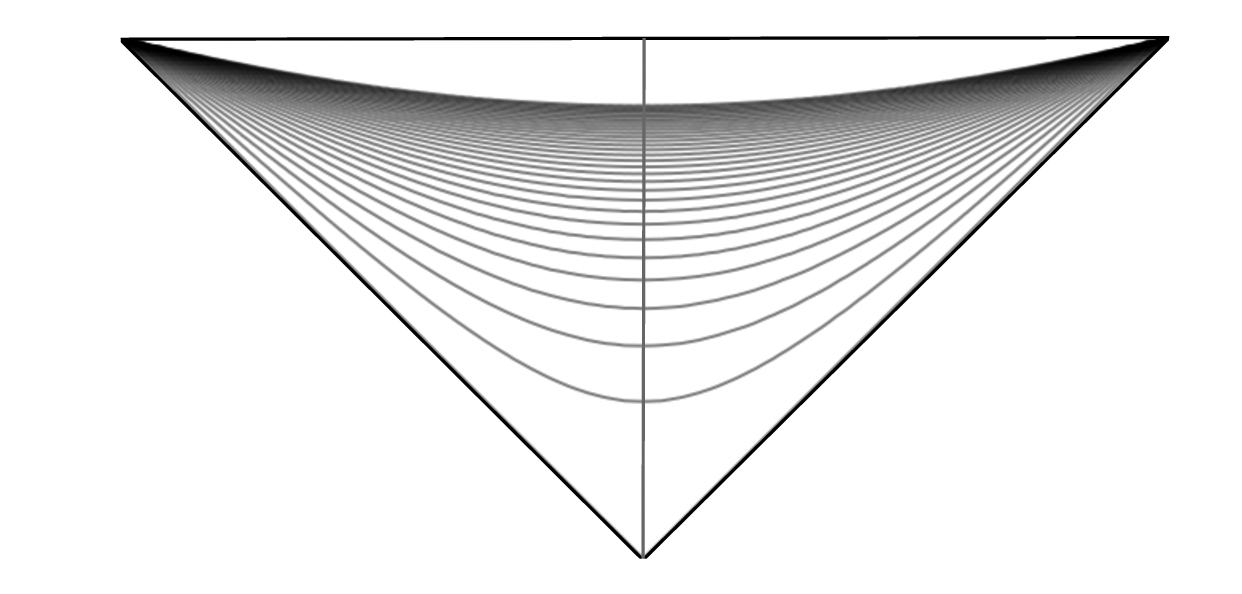}}
\newline
{\bf Figure 5.3\/} The character variety and its triple product
foliation
\end{center}

The cone point of $C$, the point $(1/2,1/2)$,
corresponds to the Fuchsian representation, the
one associated to a maximally symmetric marked
box. The diagonal edges of slope $\pm 1$, which glue together
in $C$ to make a single segment, correspond to
representations which have an extra symmetry.
The picture in Figure 1.2 comes from a point in this segment.

We call two Pappus modular groups with the same
triple invariant {\it iso-prismatic\/}, for reasons that
will become clear below.  The iso-prismatic representations
fall within the same level set.
Within each iso-prismatic family there is a unique
representation corresponding to the point on the
vertical midline in Figure 5.3.  We call this the
{\it axial representation\/}.
I might conjecture that quantities like e.g.
the Hausdorff dimension of the limit set vary
monotonically as one moves away from the
axial representation along an iso-prismatic family.

\subsection{The Geometry of Prisms}

Each marked box $M$ in the orbit
$\cal M$ has a {\it foliated flat\/} $f_M$ associated
to it.  The flat is the unique one that contains
the geodesic $\gamma_M$.  The foliation is
by geodesics parallel to $\gamma_M$.  This is the
same picture we developed in the previous chapter.
The flag $f_M$ contains the order $2$ fixed
point $p_M$ of the polarity that swaps $M$ and $i(M)$,
and $\gamma_M$ is the geodesic through $p_M$.

As in the introduction, a {\it prism\/} is a triple
of these flats corresponding to a triangle in the
Farey triangulation.  The three flats in the prism
contain a triangle of geodesics in the Farey
pattern which are pairwise one-end-asymptotic.
Now we say further that the prism is actually
foliated by such triangles of geodesics.
The way this works is that the order
$3$ isometry which, on the boundary,
permutes $i(M), t(M), b(M)$ cycles the
flats in the prism and the triangles of which
we speak are orbits of the foliating geodesics
under this action.  The triangle in the
Farey pattern is one of these triangles,
and the rest are like parallel displacements.

Our pattern of flats is really comprised of
triangle-foliated prisms.  The Farey
pattern is obtained by taking the
distinguished triangles within each
foliated prism and snapping them together.
Here is the lemma that underlies our
bending phenomenon.   By construction,
all the foliated prisms within the same
pattern are isometric to each other.
\newline
\newline
\noindent
{\bf Bending Explained:\/}
Our bending phenomenon occurs within each
iso-prismatic family.
Since the
isometry type of the prism only depends on the
triple invariant, the prisms associated to iso-prismatic
Pappus modular groups are isometric to each other.
As we vary the representations within an
iso-prismatic family, the associated order
$2$ fixed points are moving up and down the
prism. 
(Below we will see that in a certain family of
distinguished singular lines, three per prism.)
Correspondingly, the distinguished
triangles in the Farey pattern are changing
their relative positions within the prism.
Starting with one prism,  we get the adjacent ones by
reflecting isometrically across the abovementioned order $2$ points.
As the fixed points move, the geometric relation between
a prism and its neighbor changes. 
This is our bending phenomenon.  Compare
[{\bf T\/}], [{\bf P\/}], and [{\bf FG\/}].
\newline

One might wonder whether, two triangles in the foliation of
a prism are isometric.  Generically, the answer is no.
We will prove this after
a preliminary lemma.

\begin{lemma}
  \label{triple}
  Let $J$ be a triple of flags with triple invariant not equal to $\pm 1$.
  Then there
  are exactly three polarities which maps a triple of flags to itself.
  Each one permutes the
  flags by an odd permutation.
\end{lemma}

\startproof
Let $\Delta$ be the standard polarity.
We let $J'=\Delta(J).$
The vectors representing the points and
lines of $J'$ can be taken to be the same
as those representing the points and lines
of $J$.  All that has happened is that the
interpretations of these vectors changes.
If $V$ formerly represented $1$-dimensional
subspace, then $V$ is now a linear functional
representing $V^{\perp}$.

For this reason, $\chi(J')=1/\chi(J)$.
 If we permute the triple
of flags in $J'$ by an odd permutation we get a new
triple $J''$ with $\chi(J'')=1/\chi(J')=\chi(J)$.
But then there is a projective transformation $T$
such that $T(J)=J''$.
But then, $\Delta \circ T(J)=\Delta(J'')$.  This
last triple is the same as the triple $J$ but with
the flags permuted by an odd permutation.

Our duality is $\psi=\Delta \circ T$.
Note that $\psi^2(J)=J$ which means that
$\psi^2$ is a projective transformation whose
order is either $1$ or $3$.   When we
apply $\psi$ twice, we are applying the
same odd permutation, so $\psi^2$ acts
as the identity permutation on the flags 
of $J$.  Hence $\psi^2$ is the identity.
Hence $\psi$ is a polarity.
Once we have one $\psi$ we can
get two more by composing $\psi$ with
one of the order $3$ projective transformations
which cycle the flags of $J$.
\endproof

Now we  answer the question about the geometry
of the triangles within a foliated prism.

\begin{lemma}
  \label{almost}
Each  triangle in the prism
   foliation has at most one other triangle
  in the same foliation
  it could be isometric to.
  \end{lemma}

\startproof
Let $\gamma_1$ and $\gamma_2$ be
two triangles in the foliation.  Suppose
that $I$ is an isometry of $X$ such that
$I(\gamma_1)=\gamma_2$.  Since
each geodesic of $\gamma_j$ is contained
in a unique flat, we see that $I$ preserves
our triple of flats.   But then $I$ stabilizes
the triple of flags defining the triple of
flats. This forces $I$ to be one of
$6$ isometries determined by the
flag, either the identity, one of
the two projective transformations of order
$3$, or one of the polarities from Lemma \ref{triple}.
But the orbit of $\gamma_1$ under this
order $6$ group consists of either one or two
triangles.  Hence there is at most one other
choice for $\gamma_2$.
\endproof

\subsection{Axial Representations and Inflection}

The axial representations are the ones corresponding
to the equivalence class $[1/2,x)]$ in the character variety.
Again, these are the representations corresponding to
the vertical midline in Figure 5.3.

\begin{lemma}
  \label{polar}
  Suppose that $\rho$ is a polarity
  that preserves a prism. Then
  $\rho$ has a unique fixed point
  in the prism flat that it stabilizes.
\end{lemma}

\startproof
The polarity $\rho$ acts as an odd permutation
on the triple of flags
defining a prism.  The flat $F$
corresponding to the flags swapped by $\rho$ is
stabilized by $\rho$.  Also, $\rho$ reverses the medial
foliation of $F$ defined by the swapped flats.
Being an isometry, $\rho$ also
preserves the orthogonal foliation
by singular geodesics.  But since $\rho$ is
a polarity, $\rho$ must reverse the directions
of these singular geodesics.   Since $\rho$
reverses a pair of orthogonal directions on $F$,
we see that $\rho$ reverses all directions on $F$.
Hence $\rho$ has a unique fixed point on $F$.
\endproof

There is one of these fixed points for each flat
in the prism.
We call these points the {\it inflection points\/}
of the prism.

\begin{lemma}
  For an axial representation, the inflection
  points are contained in the geodesics of
  the Farey pattern and coincide with the fixed
    points of the order $2$ elements of the group.
  \end{lemma}

\startproof
We choose a markex box $M$ in the orbit
with invariant $[(x,1/2)]$.Then
$i(M)$ has box invariant $[(1/2,x)]$.
So, there is a projective
transformation which maps $M$ to $i(M)$ and
swaps the top and bottom.  Composing this
with the polarity that maps $i(M)$ to $M$, we
get a polarity $\rho$ that preserves $M$ and swaps
the top and the bottom.
The polarity $\rho$ preserves
the marked box orbit (up to switching
tops and bottoms) and therefore the
Farey pattern.
Also $\rho$ preserves
the triple of flags associated to $M$ and
thus is as in Lemma \ref{polar}.  This
forces the inflection point fixed by $\rho$ to
be in the Farey pattern.

Each geodesic $\gamma$ in the Farey pattern
also contains an order $2$ fixed point.  If
this point does not coincide with the inflection
point, then we
would have two order $2$ isometries of the Farey
pattern fixing distinct points in the same
geodesic.  The square of the
product would generate an
infinite subgroup of isometries of the
stabilizing a prism, contradicting
Lemma \ref{triple}.
\endproof

Each prism has a triple of distinguished singular
lines.  These are the lines perpendicular to the
triangles in the triangle foliation and containing
the inflection points.  We call these the {\it inflection lines\/}.

\begin{lemma}
  For any Pappus modular group whose
  triple product is not equal to $-1$, the
    order $2$ fixed points lie on the inflection lines.
  The inflection lines in the same flat, defined relative
  to adjacent prisms which contain it, coincide.
\end{lemma}

\startproof
This is a calculation.
We work with respect to the marked box $M$ shown in
Figure 3.3.  To avoid a tedious accounting of the points and
lines which I would probably do incorrectly anyway, let me
be a big vague about which matrix I am computing.

Let $\delta$ be the polarity swapping $M$ and $i(M)$.
Let $\Sigma$ denote the set of $6$ polarities associated to
the prisms associated to $M$ and $i(M)$.   Then we have
$12$ possible linear transformations of the form
$\delta \circ \sigma$ or $\sigma \circ \delta$ where $\sigma \in
\Sigma$.
One of these, which we call $T$, has the action
$$T(t)=t, \hskip 10 pt T(b)=b, \hskip 10 pt T([1:0:0])=[1:0:0], \hskip
10 pt  T([1:1:2])=[0:1:0],$$
$$[1:0:0]=\overline{su} \cap \overline{ac}, \hskip 30 pt
[1:1:2]=\overline{sa} \cap \overline{uc}, \hskip 30 pt
[0:1:0]:=\overline{sc} \cap \overline{ua}.$$
We compute that
\begin{equation}
T = \alpha \left[\matrix{\frac{-1+x+y}{x-y}& \frac{-1+3x+y-4xy}{x-y}& \frac{1-2x-y+2xy}{x-y} \cr 0 & 1 & 0 \cr 0 & 2 & -1}\right]
        \end{equation}
        Here $\alpha$ is a constant which makes $\det(T)=1$.
        The eigenvalues for $T/\alpha$ are
        $$1, \hskip 20 pt -1 \hskip 20 pt \frac{-1+x+y}{x-y}.$$
        The corresponding eigenvectors represent the points
        $t,b,[1:0:0]$.
        Notice that $T$ acts as an isometry on the flat determined by
        the top and bottom flag of $M$.  The eigensystem calculation
        identifies $T$ as a translation along the singular lines which
        limit in one direction to $[1:0:0]$.  This is only possible if
        the
        fixed point of $\delta$ and the inflection point are on the
        same
        singular line.
        \endproof

\noindent
{\bf Remark:\/}
We fix a prism $P$ and restrict our
  attention to those Pappus modular groups
  which have $P$ as part of the pattern of flats.
  This gives us instances of all the
  representations with a certain triple invariant.
  As we range over the level set, the order
  $2$ fixed points of the associated groups sweep out
  the inflection lines.
  
  We choose $3$ points in $P$ that are permuted by
    the order $3$ symmetry $g$ of $P$ and then look at
    the group generated by $g$ and
      the elliptic polarities which fix
  these $3$ points. If we choose the points to be in
  the inflection lines, we recover the Pappus modular groups
  we have been talking about.  If we pick these points
  off the inflection lines, we get something new.  Perhaps
  this  is an alternate description of the additional family of groups
  defined  in [{\bf BLV\/}].

  \subsection{Filling}
  \label{rough}

  Let me explain one way to fill in the Farey pattern
  $\Gamma_{\cal M}$ to
  get a kind of pleated surface.  We will use the
  same construction in a ``fiberwise'' way to get a kind
  of pleated $3$-manifold that fills in the pattern of flats.

  We introduce the following notation.
  \begin{itemize}
\item   Let $\tau=(\gamma_1,\gamma_2,\gamma_3)$ be a triangle in
  $\Gamma_{\cal M}$.
\item Let  $\Pi=(f_1,f_2,f_3)$ be the prism containing $\tau$.
  \item Let
    $\ell_j \subset f_j$ be the inflection line.
    \item Let $\pi_j \in \ell_j$
  be the inflection point.
\item Let $\ell$ be the singular geodesic fixed by the order $3$
  isometries of $\Pi$.
\item Let $\pi \in \ell$ be the point fixed by each of the
  $3$ order $2$ isometries of $\Pi$.
    \end{itemize}

    The following lemmas justify the last two definitions.

  \begin{lemma}
    The fixed point set of the order $3$ isometries of $\Pi$ is a
    singular geodesic $\ell$. 
    There is a unique $\pi \in \ell$ fixed by the
    order $2$ isometries of $\Pi$.
     \end{lemma}

  \startproof
  When the associated flags are as in Figure 5.1,
  the symmetry of $\Pi$ is induced by an order $3$ rotation
  about the origin.  The corresponding fixed set in $X$
  consists of standard ellipsoids with principal lengths
  $a,a,c$.  This is a singular geodesic.  The general case
  follows from symmetry.

  By symmetry each of the order $2$ symmetries
  of $\Pi$ must stabilize $\ell$.  But these symmetries
  are induced by polarities and hence reverse the
  direction of $\ell$.  Hence each one has a unique
  fixed point on $\ell$.  These fixed points all coincide
  because the isometries of $\Pi$ form a finite group.
  \endproof

  We orient the singular geodesics
  $\ell, \ell_1, \ell_2,\ell_3$ so that they
  move from $\P^*$ to $\P$.
  Let $d$ denote the signed
  distance from $\gamma_j \cap \ell_j$ to the inflection
  point $\pi_j$.    Let $\pi_d$ denote the point
   on $\ell$ whose signed distance from $\ell$ is $d$.
   Let $F_d\Pi$ be the cone of $\tau$ to $\pi_d$.
  The boundary of
  $F_d \Pi$ is $\tau$.  I think of
  $F_d\Pi$ as kind of an ideal
  triangle in $X$ but with more exotic geometry.
 $F_d\Pi$ has three analytic pieces which
    meet along medial geodesic rays emanating
    from $\pi_d$.
    
  Doing this construction for every triangle associated
  to the Farey pattern, we get our pleated surface.
  The ``ideal triangles'' in this surface fit together
  combinatorially
    just as the ideal triangles in
  the Farey triangulation fit together to give the hyperbolic plane.
  I don't  know how to prove that this surface is embedded.
  Of course, there is one case that one can completely understand:
    For the totally symmetric representation, this construction
  recreates the Farey triangulation inside an
  isometric copy $\H^2 \subset X$.

  To fill in the prism $\Pi$ we define
  $$F\Pi=\bigcup_{d \in \R} F_d\Pi.$$
  I guess that the individual $F_d \Pi$ give
  a foliation of $F\Pi$ but I don't know how to prove this.

  Making this filling for every prism associated to the
  pattern of flats, we get a kind of pleated
  $3$-manifold associated to $\Gamma_{\cal M}$.
  I don't know how to prove this thing is embedded,
  but in any case this $3$-manifold exhibits a bending
  phenomenon.  As we vary within an iso-prismatic family,
  the filled prisms do not change their geometry but pairs of adjacent
  filled prisms
  change their relationship to each other.
  \newline
  \newline
  {\bf Remark:\/}
 The  approach I take here
  is quite elementary though perhaps it is not the best approach.
  Anna Wienhard pointed out to me that cones in higher rank
  are pretty hard to control geometrically, and she
  suggested to me that perhaps the approach
  taken in [{\bf DR\/}] might work better.

  \newpage
  \section{References}

\noindent
[{\bf Bar\/}]
T. Barbot,
{\it Three dimensional Anosov Flag Manifolds\/},  Geometry $\&$ Topology (2010)
\vskip 8 pt
\noindent
[{\bf BLV\/}], T. Barbot, G. Lee, V. P. Valerio, {\it Pappus's Theorem,
  Schwartz Representations, and Anosov Representations\/},
Ann. Inst. Fourier (Grenoble) {\bf 68\/} (2018) no. 6
\vskip 8 pt
\noindent
[{\bf BCLS\/}] M. Bridgeman, D. Canary, F. Labourie, A. Samburino,
{\it The pressure metric for Anosov representations\/}, Geometry and Functional Analysis {\bf 25\/} (2015)
\vskip 8 pt
\noindent
[{\bf DR\/}] C. Davalo and J. M. Riestenberg, {\it Finite-sided
  Dirichlet domains and Anosov subgroups\/}, arXiv 2402.06408 (2024)
\vskip 8 pt
\noindent
[{\bf FG\/}]. V. Fock and A. Goncharov, {\it Moduli Spaces of local
  systems and
  higher Teichmuller Theory\/}, Publ. IHES {\bf 103\/} (2006)
\vskip 8 pt
\noindent
[{\bf G\/}] W. Goldman, {\it Convex real projective structures on compact surfaces\/}, J.
Diff. Geom. {\bf 31\/} (1990)
\vskip 8 pt
\noindent
[{\bf GW\/}] O. Guichard, A. Wienhard, {\it Anosov Representations: Domains of
  Discontinuity and applications\/}, Invent Math {\bf 190\/} (2012)
\vskip 8 pt
\noindent
[{\bf Hit\/}] N. Hitchin, {\it Lie Groups and Teichmuller Space\/}, Topology {\bf 31} (1992)
\vskip 8 pt
\noindent
[{\bf KL\/}] M. Kapovich, B. Leeb, {\it Relativizing characterizations of
  Anosov subgroups, I (with an appendix by Gregory A. Soifer).\/} Groups Geom. Dyn.  {\bf 17\/}  (2023)
  \vskip 8 pt
  \noindent
[{\bf Lab\/}] 
F. Labourie, {\it Anosov Flows, Surface Groups and Curves in Projective Spaces\/},
P.A.M.Q {\bf 3\/} (2007)
\vskip 8 pt
\noindent
[{\bf P\/}] R. Penner, {\it The Decorated Teichmuller Theory of
  punctured surfaces\/}, Comm. Math. Pys. {\bf 113\/} (1987)
\vskip 8 pt
\noindent
[{\bf S0\/}] R. E. Schwartz, {\it Pappus's Theorem and the Modular Group\/},
Publ. IHES (1993)
\vskip 8 pt
\noindent
[{\bf T\/}] W. Thurston, {\it The Geometry and Topology of Three
  Manifolds\/}, Princeton University Notes (1978)

\end{document}